\documentclass[preprint,12pt]{article}

\usepackage{arxiv}
\usepackage[english]{babel}
\usepackage[utf8x]{inputenc}
\usepackage[T1]{fontenc}
\usepackage{amsmath,amsfonts,amssymb,systeme} 
\usepackage{graphicx} 
\usepackage{lineno} 
\usepackage{epstopdf} 
\usepackage{subcaption} 
\usepackage{hyperref}
\usepackage{cleveref}
\usepackage{xcolor}
\usepackage{makecell}

\DeclareMathOperator{\w}{\omega}
\DeclareMathOperator{\wnot}{\omega_0}
\DeclareMathOperator{\wk}{\omega_k}
\DeclareMathOperator{\W}{\Omega}
\DeclareMathOperator{\atantwo}{atan2}

\graphicspath{{img/}}

\title{Resonant phase lags of a Duffing oscillator}
\author{Martin Volvert\\
	Space Structures and Systems Laboratory,\\
	Aerospace and Mechanical Engineering Department,\\
	University of Liège, Belgium\\
	\texttt{m.volvert@uliege.be}
	\and Gaëtan Kerschen\\
	Space Structures and Systems Laboratory,\\
	Aerospace and Mechanical Engineering Department,\\
	University of Liège, Belgium\\
	\texttt{g.kerschen@uliege.be}}
\begin{document}
    \maketitle
    \section{Introduction}
\label{sec:1_Intro}

The resonant behavior of linear systems can be characterized either with the concept of an amplitude resonance or a phase resonance. Amplitude resonance corresponds to a relative maximum in the frequency response function whereas phase resonance is associated with quadrature between the displacement and the external forcing. At phase resonance, the external forcing cancels exactly the damping force with the result that the resonance frequency coincides with the natural frequency of the linear system. The difference between the two resonances remains small for weakly damped systems. Phase resonance-based testing \cite{WRIGHT} which excites the individual modes of the system in turn was largely exploited during the early days of experimental modal analysis because it provides accurate estimation of the modal parameters. With the advent of advanced system identification techniques such as the stochastic subspace identification method \cite{DEMOOR}, phase resonance testing has been less and less employed for linear modal analysis. 

For nonlinear systems, the phase lag quadrature criterion was first extended to synchronous motions using harmonic balance in \cite{JSVPEETERS} and then to arbitrary periodic motions using Melnikov analysis in \cite{HALLER}. These efforts triggered the development of nonlinear phase resonance testing which targets the identification of the nonlinear normal modes (NNMs) defined as periodic solutions of the unforced, undamped system \cite{VAKAKISBOOK,MSSPKERSCHEN}. Basic \cite{MSSPPEETERS,NEILD,ALLEN} and more advanced (control-based) strategies \cite{RENSON,PETER,PETER2,THOMAS,THOMAS2,OZGUVEN,YABUNO,SCHEEL,ABELOOS} were developed during the last decade. In this context, phase-locked loops (PLLs) are particularly effective for tracking phase quadrature for increasing forcing amplitudes, as first proposed in \cite{PETER}. In addition, like control-based continuation \cite{RENSON}, phase control may also stabilize unstable periodic solutions.  

Despite the great promise of PLLs for experimental modal analysis of nonlinear systems, two difficulties remain for an accurate and thorough characterization of nonlinear resonant behaviors. First, according to \cite{JSVPEETERS}, the correspondence between the quadrature curves identified using PLLs and NNMs is only valid for multi-harmonic forcing. In the presence of modal interactions, the discrepancy can be very important \cite{RENSON2}. Second, nonlinear systems can exhibit additional resonances including superharmonic and subharmonic resonances. Even if recent theoretical \cite{HALLER} and numerical \cite{JSVVOLVERT} studies investigated these secondary resonances under the banner of nonlinear modes, it is not yet fully clear how they can be identified using phase resonance testing.

To provide a solid theoretical framework for the use of PLLs in nonlinear experimental modal analysis, the present study aims to pursue the study initiated in \cite{Leung}. Specifically, we revisit the resonant behavior of a harmonically-forced Duffing oscillator with a specific attention to phase resonance and to its relation with amplitude resonance. To this end, the different families of resonances including primary (1:1), superharmonic (k:1), subharmonic (1:$\nu$) and ultra-subharmonic (k:$\nu$) resonances are carefully studied using first and higher-order averaging. 

The paper is organized as follows. Section 2 briefly recalls the principles behind averaging in nonlinear dynamics. Section 3 focuses on the amplitude and phase resonances of both a linear and a Duffing oscillator whereas Section 4 extends the investigations to different secondary resonances. In Section 5, the findings obtained through the analytical derivations are verified using numerical simulations. The conclusions of the present study are summarized in Section 6.

    \section{General Approach}

\subsection{Resonances of a Duffing oscillator}

The governing equation of motion of a harmonically-forced Duffing oscillator is
\begin{equation}
    m\ddot{x}(t)+c\dot{x}(t)+kx(t)+k_{nl}x^3(t)=f\sin{\omega t}
    \label{eq:DUFFING}
\end{equation}
where $m$, $c$, $k$ and $k_{nl}$ represent the mass, damping, linear and nonlinear stiffness coefficients, respectively. $f$ is the forcing amplitude whereas $\omega$ is the excitation frequency of period $T$. The natural frequency of the undamped, linearized system is $\omega_0=\sqrt{\frac{k}{m}}$ \cite{KOVACIC}. The Duffing oscillator is said to be hardening when $k_{nl}>0$, and softening when $k_{nl}<0$. Only the hardening case is studied here.

Through mass normalization, Equation \eqref{eq:DUFFING} can be recast into:

\begin{equation}
    \ddot{x}(t)+2\bar{\zeta}\wnot\dot{x}(t)+\omega_0^2x(t)+\alpha x^3(t)=\bar{\gamma}\sin{\omega t}
    \label{eq:DUFFING_NORM}
\end{equation}
where $\bar{\zeta}=\frac{c}{2\sqrt{km}}$, $\alpha=k_{nl}/m$ and $\bar{\gamma}=f/m$. The coefficients are set to $m=1$kg, $c=0.01$kg/s, $k=1$N/m and $k_{nl}=1$N/m$^3$ throughout the present study.

If we consider the Fourier decomposition of the displacement
\begin{equation}
x(t) = A_0 + \sum_{k=1}^n A_k \sin \left(\omega_k t-\phi_k\right)
\label{eq:FourierDispPhase}
\end{equation}
where $\omega_k=\frac{k\omega}{\nu}$ (with $\nu$ a positive integer), $A_k$ and $\phi_k$ are the frequency, amplitude and phase lag  of the $k$-th harmonic of the displacement, respectively, then Equation (\ref{eq:FourierDispPhase}) shows that each harmonic $k$ may trigger a resonance if $\omega_k$ corresponds to the (amplitude-dependent) frequency of the primary resonance of the system. These resonances can be divided into four categories, namely $1:1$ primary resonance ($k=\nu=1$), $k:1$ superharmonic resonances, $1:\nu$ subharmonic resonances and $k:\nu$ ultra-subharmonic resonances. 

\subsection{Averaging around the \texorpdfstring{$k:\nu$}{} resonance}

We consider a weakly nonlinear oscillator of the form:

\begin{equation}
    \ddot{x}(t)+\omega_0^2 x(t) = \varepsilon f(x(t),\dot{x}(t))
    \label{eq:2ndeq}
\end{equation}
When $\varepsilon=0$, the periodic solution of \eqref{eq:2ndeq} is written as:

\begin{equation}
    x(t) = u\cos{\wnot t} - v\sin{\wnot t}
    \label{eq:xt}
\end{equation}
where $u$ and $v$ are constants. When $\varepsilon\neq 0$, we seek a solution of frequency $\wk$ such that $\wk^2-\wnot^2=\varepsilon\W$. The solution is expressed as in Equation \eqref{eq:xt} but with time-dependent $u$ and $v$:
\begin{equation}
    x(t) = u(t)\cos{\wk t} - v(t)\sin{\wk t}
    \label{eq:xt2}
\end{equation}
We impose that the velocity should have the same form as in the case $\varepsilon=0$, \textit{i.e.},
\begin{equation}
    \dot{x}(t) = -u(t)\wk\sin{\wk t} - v(t)\wk\cos{\wk t}
    \label{eq:dxt}
\end{equation}
Equation \eqref{eq:dxt} holds if:
\begin{equation}
    \dot{u}(t)\cos{\wk t}-\dot{v}(t)\sin{\wk t} = 0
    \label{eq:CONDVDP}
\end{equation}
Differentiating Equation \eqref{eq:dxt} and replacing $\ddot{x}(t)$ and $x(t)$ in Equation \eqref{eq:2ndeq} yields:
\begin{equation}
    \dot{u}(t)\wk\sin{\wk t} + \dot{v}(t)\wk\cos{\wk t} = -\varepsilon \left[f(x(t),\dot{x}(t))+\wk\W x(t)\right]
    \label{eq:2ndequv}
\end{equation}
Finally, taking into account Equations \eqref{eq:CONDVDP} and \eqref{eq:2ndequv} and solving for $\dot{u}$ and $\dot{v}$, a system of first-order equations is obtained:
\begin{equation}
    \begin{cases}
    \dot{u}=-\frac{\varepsilon}{\wk}\left[f(x(t),\dot{x}(t))+\wk\W x(t)\right]\sin{\wk t}\\
    \dot{v}=-\frac{\varepsilon}{\wk}\left[f(x(t),\dot{x}(t))+\wk\W x(t)\right]\cos{\wk t}
    \end{cases}
    \label{eq:1steq}
\end{equation}
This system has a suitable form to apply first- or higher-order averaging. First-order averaging is performed herein using the Krylov-Bogolyubov technique \cite{Krylov, NAYFEH2}. Higher-order averaging is based on the Lie transform algorithm \cite{NAYFEH2}; it was implemented by Yagasaki in the \textit{haverage.m} Mathematica package \cite{YAGA1,YAGA2}. 

$x(t)$ is often represented using the polar coordinates $r$ and $\phi$ such that $x(t) = r(t)\sin{(\w t-\phi(t))}$ with $r=\sqrt{u^2+v^2}$ and $\phi=\atantwo(-u,-v)$. For conciseness, the time dependence for $u$, $v$, $r$ and $\phi$ is dropped in the remainder of this article.
    \section{Primary Resonance (\texorpdfstring{$k=\nu=1$}{})}

Considering Equation \eqref{eq:DUFFING_NORM}, we scale the system such that $\bar{\zeta}=\varepsilon\zeta$ and $\bar{\gamma}=\varepsilon^{3/2}\gamma$, with $\zeta$, $\gamma=\mathcal{O}(1)$. If $x=\sqrt{\varepsilon}y$, we obtain:
\begin{equation}
    \ddot{y}(t)+2\varepsilon\zeta\wnot\dot{y}(t)+\omega_0^2 y(t)+\varepsilon\alpha y^3(t)=\varepsilon\gamma\sin{\omega t}
\end{equation}
The forcing frequency is in the vicinity of the natural frequency of the linear system, i.e., $\w^2-\wnot^2=\varepsilon\Omega$. The displacement is expressed as:
\begin{equation}
    x(t) = \sqrt{\varepsilon}r\sin{(\w t-\phi)} = A\sin{(\w t-\phi)}
\end{equation}

\subsection{Linear system}
Applying first-order averaging to the linear system ($\alpha=0$) yields:
\begin{equation}
    \begin{cases}
    \dot{r}=-\frac{\varepsilon}{\w}\left(2\zeta\wnot\w r-\gamma\sin{\phi}\right)\\
    \dot{\phi}=\frac{\varepsilon}{\w}\left(\W+\gamma\cos{\phi}\right)
    \end{cases}
    \label{eq:drdphiPrimLin}
\end{equation}
Assuming a steady-state response, \textit{i.e.}, $\dot{r}=\dot{\phi}=0$, the motion around resonance is governed by:
\begin{equation}
    \begin{cases}
    2\zeta\wnot\w r=\gamma\sin{\phi}\\
    -\W r = \gamma\cos{\phi}
    \end{cases}
    \label{eq:PrimLin}
\end{equation}

The resonant behavior of a linear oscillator can be described in two different ways, i.e., either when the amplitude of the frequency response undergoes a relative maximum (i.e., amplitude resonance denoted by a subscript \textit{a} herein) or when the displacement is in quadrature with the external forcing (i.e., phase resonance denoted by a subscript \textit{p}). Both cases are detailed in what follows.

\subsubsection{Phase lag at amplitude resonance}
Amplitude resonance occurs when both $\frac{\partial r}{\partial\w}$ and $\frac{\partial r}{\partial\phi}$ are equal to 0. From Equation \eqref{eq:PrimLin}, we obtain:

\begin{equation}
    \begin{cases}
    \frac{\partial r}{\partial \phi} = \frac{\gamma}{2\zeta\wnot\w}\left(\cos{\phi}-\frac{\sin{\phi}}{\w} \frac{\partial\w}{\partial\phi}\right)=0\\
    \frac{\partial r}{\partial \w} =\frac{\gamma}{2\zeta\wnot\w}\left(\cos{\phi}\frac{\partial\phi}{\partial\w}-\frac{\sin{\phi}}{\w}\right)=0
    \end{cases}
    \label{eq:drdphi_drdw_prim}
\end{equation}
Both relations are equivalent. The second relation of Equations \eqref{eq:PrimLin} provides an expression for $\omega$:
\begin{equation}
    \w=\sqrt{\wnot^2-\frac{\varepsilon\gamma}{r}\cos{\phi}}
\end{equation}
$\frac{\partial \w}{\partial \phi}$ is obtained by isolating $\W$ in the second equation of \eqref{eq:PrimLin} and making use of the chain rule $\frac{\partial\w}{\partial\phi} = \frac{\partial\w}{\partial\W}\frac{\partial\W}{\partial\phi}$ such that
\begin{equation}
    \frac{\partial\w}{\partial\phi} = \frac{\varepsilon}{2\w}\left(\frac{\gamma}{r}\sin{\phi}+\frac{\gamma}{r^2}\cos{\phi}\frac{\partial r}{\partial \phi}\right)
\end{equation}
which can be inserted in the first relation of \eqref{eq:drdphi_drdw_prim}:
\begin{equation}
    \frac{\partial r}{\partial \phi} = \frac{\gamma\sin{\phi}}{2\zeta\wnot\w}\frac{\left(\w-\varepsilon\zeta\wnot\tan{\phi}\right)}{\left(\w\tan{\phi}+\varepsilon\zeta\wnot\right)}=0
    \label{eq:drdphi_PrimLin}
\end{equation}
This relationship is satisfied when the phase lag  takes the form:
\begin{equation}
    \tan{\phi_a}=\frac{\w_a}{\varepsilon\zeta\wnot}=\frac{\sqrt{1-2\bar{\zeta}^2}}{\bar{\zeta}}
\end{equation}
The corresponding frequency and amplitude are
\begin{equation}
    \w_a=\wnot\sqrt{1-2\bar{\zeta}^2},\;\; A_a = \frac{\bar{\gamma}}{2\bar{\zeta}\wnot^2\sqrt{1-\bar{\zeta}^2}}
    \label{eq:phi_linear}
\end{equation}
It should be noted that $\phi=0$ or $\phi=\pi$ also verify Equation \eqref{eq:drdphi_PrimLin} and correspond to the purely static and inertial responses, respectively.
\subsubsection{Phase lag at phase resonance}
Phase quadrature $\phi_p=\frac{\pi}{2}$ occurs when the excitation frequency corresponds to the natural frequency of the undamped system, i.e., when $\w=\w_p=\wnot$. The amplitude at phase resonance is $A_p=\frac{\bar{\gamma}}{2\bar{\zeta}\wnot^2}$.

\subsection{Nonlinear system}
First-order averaging applied to the nonlinear system ($\alpha \neq 0$) gives:
\begin{equation}
    \begin{cases}
    \dot{r}=-\frac{\varepsilon}{\w}\left(\zeta\wnot\w r-\frac{\gamma}{2}\sin{\phi}\right)\\
    \dot{\phi}=-\frac{\varepsilon}{\w}\left(\frac{\alpha}{8}\left(3r^2-\frac{4\W}{\alpha}\right)-\frac{\gamma}{2}\cos{\phi}\right)
    \end{cases}
    \label{eq:drdphiPrim}
\end{equation}
The steady-state solutions around the primary resonance are governed by:
\begin{equation}
    \begin{cases}
    \zeta\wnot \w r=\frac{\gamma}{2}\sin{\phi}\\
    \frac{\alpha}{8}\left(3r^2-\frac{4\W}{\alpha}\right)r = \frac{\gamma}{2}\cos{\phi}
    \end{cases}
    \label{eq:drdphiPrim0}
\end{equation}

\subsubsection{Phase lag at amplitude resonance}
Following the same procedure as for the linear system, we obtain:
\begin{equation}
    \begin{cases}
    \frac{\partial r}{\partial \phi} = \frac{\gamma}{2\zeta\wnot\w}\left(\cos{\phi}-\frac{\sin{\phi}}{\w} \frac{\partial\w}{\partial\phi}\right)=0\\
    \frac{\partial r}{\partial \w} =\frac{\gamma}{2\zeta\wnot\w}\left(\cos{\phi}\frac{\partial\phi}{\partial\w}-\frac{\sin{\phi}}{\w}\right)=0
    \end{cases}
    \label{eq:drdphidrdwPrim}
\end{equation}
and
\begin{equation}
    \frac{\partial\w}{\partial\phi}=\frac{\varepsilon}{2\w}\left(\left[\frac{6\alpha r}{4} + \frac{\gamma}{r^2}\cos{\phi}\right]\frac{\partial r}{\partial \phi} + \frac{\gamma}{r}\sin{\phi}\right)
\end{equation}
Eventually,
\begin{equation}
    \frac{\partial r}{\partial \phi} = \frac{4\zeta\wnot\w^2\gamma\sin{\phi}(\w-\varepsilon\zeta\wnot\tan{\phi})}{8\zeta^2\wnot^2\w^4\tan{\phi}+\varepsilon(3\alpha\gamma^2\sin{\phi}^2\tan{\phi}+8\zeta^3\wnot^3\w^3)}
\end{equation}
This relation is verified when:
\begin{equation}
    \tan{\phi_a}=\frac{\w_a}{\varepsilon\zeta\wnot}
    \label{eq:Prim_tanPhi}
\end{equation}
and, from \eqref{eq:drdphiPrim0} and \eqref{eq:Prim_tanPhi}, it is possible to derive $A_a$, $\w_a$ and $\phi_a$ as a function of the forcing and the system parameters:
\begin{equation}
    \begin{cases}
        A_a = \sqrt{\frac{2\wnot^2}{3\alpha}
              \left((\bar{\zeta}^2-1) + \sqrt{(1-\bar{\zeta}^2)^2+\frac{3\alpha \bar{\gamma}^2}{4\bar{\zeta}^2\wnot^6}}\right)}\\
        \w_a = \frac{\wnot}{\sqrt{2}}\sqrt{1-3\bar{\zeta}^2+\sqrt{\left(1-\bar{\zeta}^2\right)^2+\frac{3\alpha\bar{\gamma}^2}{4\bar{\zeta}^2\wnot^6}}}\\
        \tan{\phi_a}=\frac{\sqrt{1-3\bar{\zeta}^2+\sqrt{\left(1-\bar{\zeta}^2\right)^2+\frac{3\alpha\bar{\gamma}^2}{4\bar{\zeta}^2\wnot^6}}}}{\sqrt{2}\bar{\zeta}}.
    \end{cases}
    \label{eq:Prim_AmpRes}
\end{equation}

\subsubsection{Phase lag at phase resonance}
Imposing $\phi_p=\pi/2$ in Equations \eqref{eq:drdphiPrim0} yields:
\begin{equation}
    \begin{cases}
    A_p=\frac{\bar{\gamma}}{2\bar{\zeta}\wnot\w_p}\\
    \w_p= \wnot\sqrt{1+\frac{3\alpha}{4\wnot^2}A^2_p}
    \end{cases}
    \label{eq:PRNMPrim}
\end{equation}
from which the expressions of the amplitude and frequency at phase resonance can be deduced:
\begin{equation}
    A_p = \sqrt{\frac{2\wnot^2}{3\alpha}\left(\sqrt{1+\frac{3\alpha\bar{\gamma}^2}{4\bar{\zeta}^2\wnot^6}}-1\right)}
    \label{eq:PRIM_rp}
\end{equation}
and
\begin{equation}
    \w_p = \frac{\wnot}{\sqrt{2}}\sqrt{1+\sqrt{1+\frac{3\alpha\bar{\gamma}^2}{4\bar{\zeta}^2\wnot^6}}}
    \label{eq:PRIM_wp}
\end{equation}
We note that Equations (\ref{eq:PRNMPrim}) correspond to those that would be obtained by applying the energy balance principle \cite{HILL,SALLES} to the NNMs of the undamped, unforced system and neglecting higher-order harmonics. Under this latter assumption, this means that phase resonance testing amounts to exciting the underlying NNMs.

\subsubsection{Discussion}

This section has derived analytical expressions of the amplitude, frequency and phase of a Duffing oscillator at amplitude and phase resonances. Of specific interest is the difference in frequency between amplitude and phase resonances, $\Delta\omega=\w_p-\w_a$:
\begin{equation}
    \Delta \w = \frac{\wnot}{\sqrt{2}}\left(\sqrt{1+\sqrt{1+\frac{3\alpha\bar{\gamma}^2}{4\bar{\zeta}^2\wnot^6}}}-\sqrt{1-3\bar{\zeta}^2+\sqrt{\left(1-\bar{\zeta}^2\right)^2+\frac{3\alpha\bar{\gamma}^2}{4\bar{\zeta}^2\wnot^6}}}\right)\label{DW}
\end{equation}
We note that the phase resonance of a harmonically-forced Duffing oscillator is rarely discussed in the technical literature. The reason might come from the fact that perturbation techniques do not always make a distinction between amplitude ans phase resonances. For example, the method of multiple scales \cite{NAYFEH} yields around the primary resonance:
\begin{equation}
    \begin{cases}
    \zeta \wnot^2 r=\frac{\gamma}{2}\sin{\phi}\\
    \frac{\alpha}{8}\left(3r^2-\frac{8\wnot(\w-\wnot)}{\varepsilon\alpha}\right)r = \frac{\gamma}{2}\cos{\phi}
    \end{cases}
    \label{eq:Nayfeh}
\end{equation}
Since the amplitude $r$ is maximum when $\phi=\pi/2$, amplitude and phase resonances are predicted to occur simultaneously with: 
\begin{equation}
    \begin{cases}
    A=\frac{\bar{\gamma}}{2\bar{\zeta}\wnot^2}\\
    \w=\wnot+\frac{3\alpha\bar{\gamma}^2}{32\bar{\zeta}^2\wnot^5}
    \end{cases}
    \label{eq:Nayfeh_AmpRes}
\end{equation}
Interestingly, the multiple scales method predicts that the amplitude at resonance of the Duffing oscillator is identical to that of the phase resonance of the underlying linear system.

\begin{figure}[htpb] 
  \begin{subfigure}[b]{0.5\linewidth}
    \centering
    \includegraphics[width=1\linewidth]{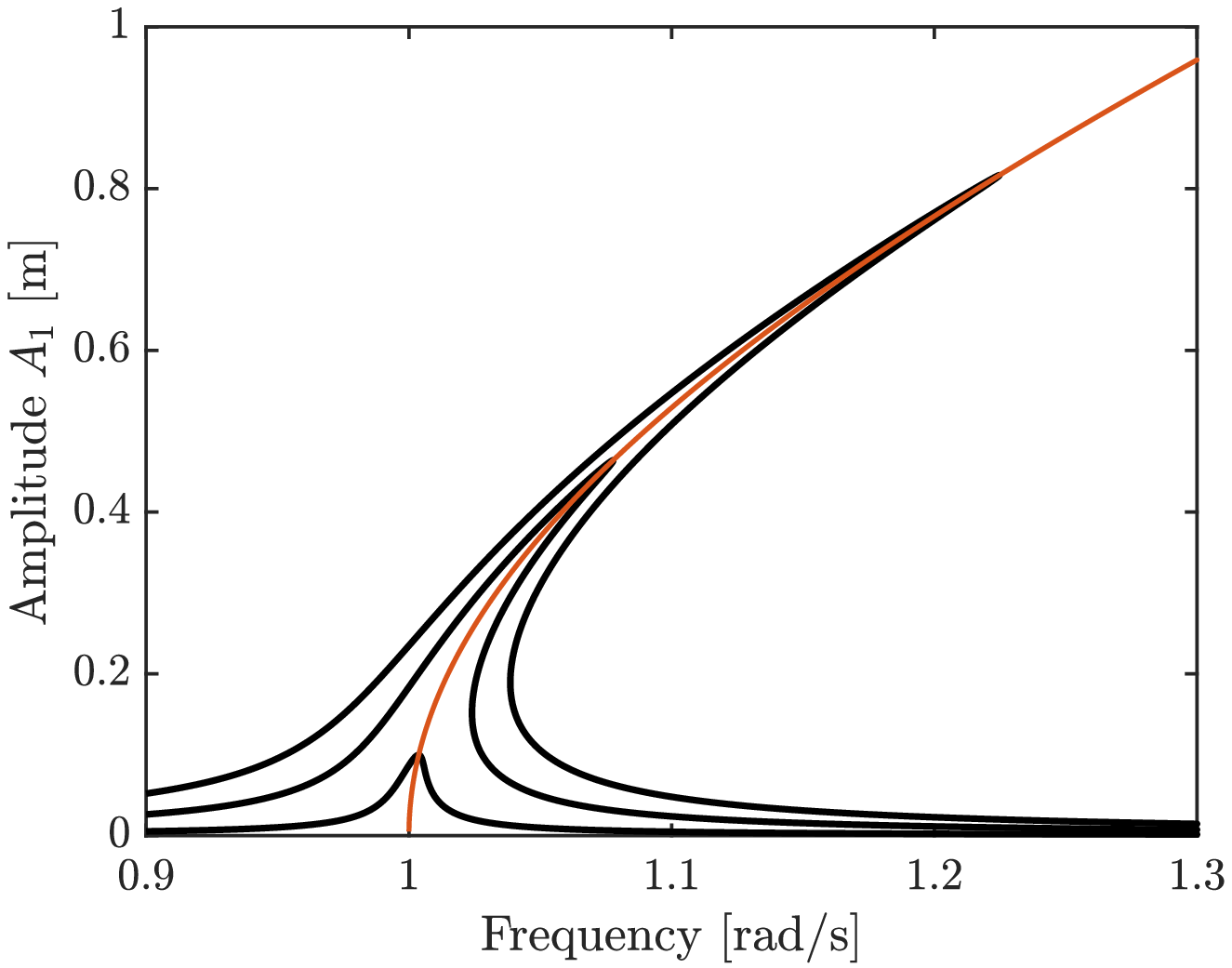} 
    \caption{\label{fig:1_NFRC}}
  \end{subfigure}
  \begin{subfigure}[b]{0.5\linewidth}
    \centering
    \includegraphics[width=1\linewidth]{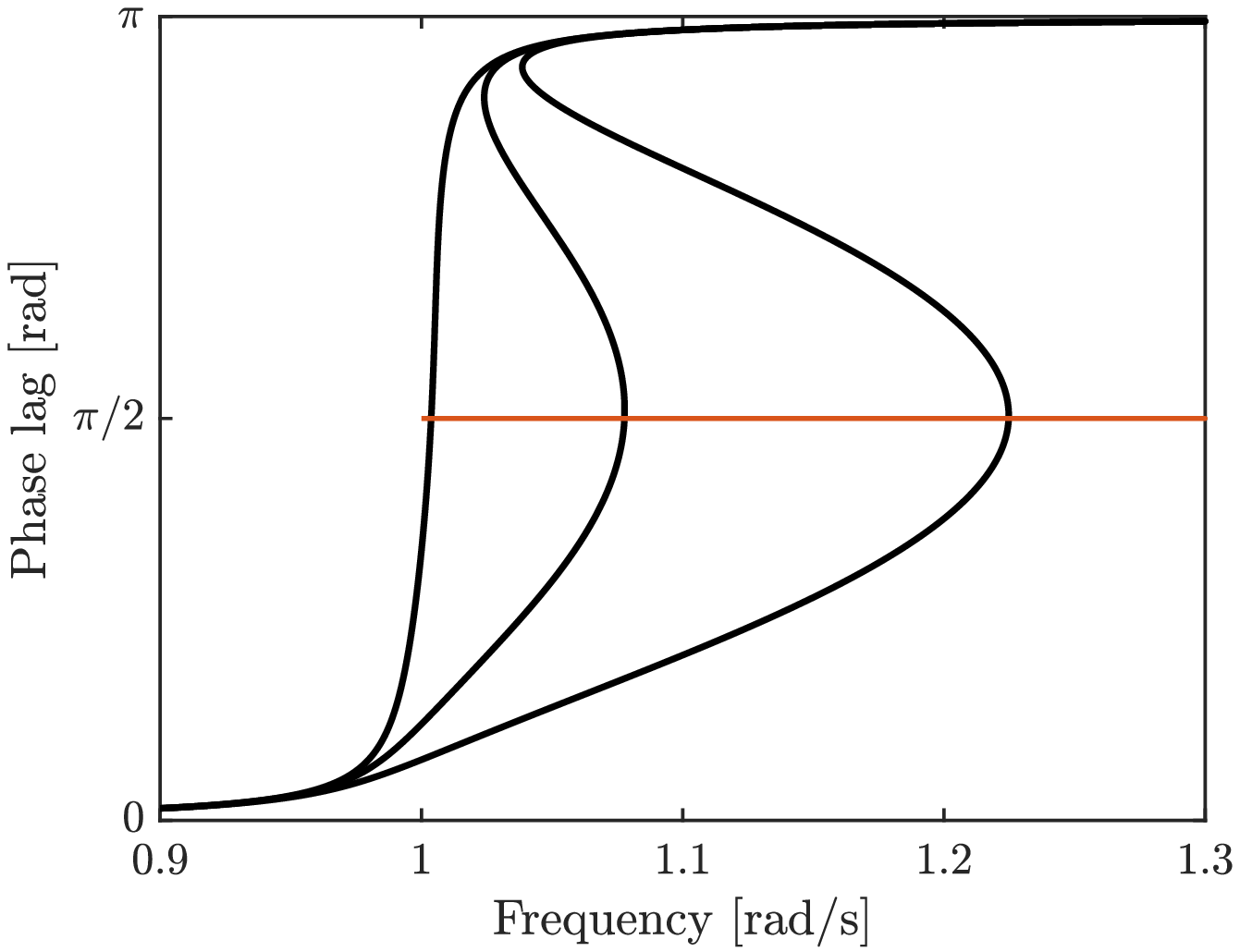}
    \caption{\label{fig:1_PHASE}}
  \end{subfigure} 
  \caption{Nonlinear frequency responses (black) around the primary resonance of the Duffing oscillator for forcing amplitudes of 0.001N, 0.005N and 0.01N and the phase resonance curve (orange): (\subref{fig:1_NFRC}) amplitude and (\subref{fig:1_PHASE}) phase lag.}
  \label{fig:PRIMRES_H1}
\end{figure}

Getting back to Equation (\ref{DW}) and performing a Taylor series expansion indicates that the frequency difference is of the order of $\mathcal{O}\left(\bar{\zeta}^2\right)$, as in the linear case. For weak to moderate damping (i.e., not beyond a few percent, which is usually the case for mechanical structures), it thus follows that phase resonance lies in the immediate neighborhood of amplitude resonance. This is illustrated in Figures \ref{fig:1_NFRC} and \ref{fig:1_PHASE}, where the phase resonance curve constructed thanks to Equations \eqref{eq:PRIM_rp} and \eqref{eq:PRIM_wp} is superposed to the nonlinear frequency responses calculated from Equations \eqref{eq:drdphiPrim0} for different forcing amplitudes. Expression (\ref{DW}) also evidences that a hardening nonlinearity ($\alpha$>0) brings amplitude and phase resonances closer to each other, and conversely for a softening nonlinearity.

    \section{Secondary Resonances}
Considering the mass-normalized equation of the Duffing oscillator \eqref{eq:DUFFING_NORM}, we scale the system such that $\bar{\zeta}=\varepsilon^d\zeta$ and $\bar{\gamma}=\sqrt{\varepsilon}\gamma$, with $\mu$, $\gamma=\mathcal{O}(1)$, and $d$ is a positive integer. If we let $x=\sqrt{\varepsilon}y$, we obtain
\begin{equation}
    \ddot{y}(t)+2\varepsilon^d\zeta\wnot\dot{y}(t)+\omega_0^2 y(t)+\varepsilon\alpha y^3(t)=\gamma\sin{\omega t}
    \label{eq:EOM1D_SCALE}
\end{equation}
If $\varepsilon=0$, then Equation \eqref{eq:EOM1D_SCALE} has a periodic solution $y(t)=\Gamma \sin{\omega t}$ with $\Gamma=\gamma/(\omega_0^2-\omega^2)$. Introducing $z(t)=y(t)-\Gamma \sin{\omega t}$ in Equation \eqref{eq:EOM1D_SCALE} yields a weakly nonlinear oscillator with a form suitable for first- or higher-order averaging:

\begin{equation}
    \ddot{z}(t) + \omega_0^2 z(t) = \varepsilon f(z(t),\dot{z}(t),\omega t, \varepsilon)
    \label{eq:WEAKNONLINEQ}
\end{equation}
where
\begin{equation}
    f(z(t),\dot{z}(t),\omega t, \varepsilon) = -\alpha(z(t)+\Gamma\sin{\omega t})^3 - 2\varepsilon^{d-1}\zeta\wnot(\dot{z}(t)+\omega\Gamma\cos{\omega t})
    \label{eq:NONLINTERM}
\end{equation}
where the forcing frequency is close to a fraction of the natural frequency of the linear system, i.e., $\wk^2-\wnot^2=\varepsilon\W$. The solution around the $k:\nu$ resonance can therefore be expressed as:
\begin{equation}
    x(t) = \sqrt{\varepsilon}\left(r\sin{(\wk t-\phi)}+\Gamma\sin{\w t}\right) = A_k\sin{(\wk t-\phi)}+\bar{\Gamma}\sin{\w t}
\end{equation}
where $A_k=\sqrt{\varepsilon}r$ and $\bar{\Gamma}=\bar{\gamma}/(\omega_0^2-\omega^2)$.

\subsection{\texorpdfstring{$3:1$}{} resonance}
The first secondary resonance studied is the $3:1$ superharmonic resonance, \textit{i.e.}, $k=3$ and $\nu=1$. Using Equation \eqref{eq:EOM1D_SCALE} with $d=1$, first-order averaging provides:
\begin{equation}
    \begin{cases}
    \dot{r}=-\varepsilon\left(\zeta\wnot r - \frac{\alpha\Gamma^3}{24\w}\sin{\phi}\right)\\
    \dot{\phi}=-\varepsilon\left(\frac{\alpha}{24\w}\left(3r^2+6\Gamma^2-\frac{4\W}{\alpha}\right)-\frac{\alpha\Gamma^3}{24\w r}\cos{\phi}\right)
    \end{cases}
    \label{eq:drdphi31}
\end{equation}
Assuming a steady-state solution, $\dot{r}=\dot{\phi}=0$, we have:
\begin{equation}
    \begin{cases}
    \frac{24\zeta\wnot \w }{\alpha}r =\Gamma^3\sin{\phi}\\
    \left(3r^2+6\Gamma^2-\frac{4\W}{\alpha}\right)r=\Gamma^3\cos{\phi}
    \end{cases}
    \label{eq:drdphi31_0}
\end{equation}

\subsubsection{Phase lag at amplitude resonance}
As for the primary resonance, amplitude resonance for the $3:1$ resonance occurs when $\frac{\partial r}{\partial \w} = \frac{\partial r}{\partial \phi} = 0$:

\begin{equation}
    \begin{cases}
    \frac{\partial r}{\partial \phi} = \frac{\alpha\Gamma^3}{24\zeta\wnot\w}\left(\left[\frac{7\w^2-\wnot^2}{\w(\wnot^2-\w^2)}\right]\sin{\phi}\frac{\partial\w}{\partial\phi}+\cos{\phi}\right)=0\\
    \frac{\partial r}{\partial \w} =\frac{\alpha\Gamma^3}{24\zeta\wnot\w}\left(\left[\frac{7\w^2-\wnot^2}{\w(\wnot^2-\w^2)}\right]\sin{\phi}+\cos{\phi}\frac{\partial\phi}{\partial\w}\right)=0
    \end{cases}
    \label{eq:drdphi_drdw_31}
\end{equation}

Isolating $\W$ from \eqref{eq:drdphi31_0}, using the chain rule $\frac{\partial\w}{\partial\phi} = \frac{\partial\w}{\partial\W}\frac{\partial\W}{\partial\phi}$ and inserting it in $\frac{\partial r}{\partial \phi}$  gives:
\begin{equation}
    \frac{\partial r}{\partial \phi} = \frac{\alpha\Gamma^3\left(\frac{(7\w^2-\wnot^2)}{\w(\wnot^2-\w^2)}\frac{\varepsilon\zeta\wnot}{3}\sin{\phi} + \left(1-\varepsilon\frac{\alpha\Gamma^2}{3(\wnot^2-\w^2)}+\varepsilon\frac{2\zeta\wnot\w}{\wnot^2-\w^2}\frac{1}{\tan{\phi}}\right)\cos{\phi}\right)}{24\zeta\wnot\w\left(1-\varepsilon\left(\frac{\alpha\Gamma^2}{3(\wnot^2-\w^2)}+\frac{2\zeta\wnot\w}{\wnot^2-\w^2}\frac{1}{\tan{\phi}}-\frac{\alpha^3\Gamma^3\sin{\phi}^2}{6912\zeta^2\wnot^2\w^3}-\frac{\zeta\wnot}{3\sin{\phi}\tan{\phi}}\right)\right)}
\end{equation}
The numerator is 0 when:
\begin{equation}
    \tan^2{\phi} - \frac{3\w(\wnot^2-\w^2)}{\varepsilon\zeta\wnot(\wnot^2-7\w^2)} \left(1-\frac{\varepsilon\alpha\Gamma^2}{3(\wnot^2-\w^2)}\right)\tan{\phi}-\frac{6\w^2}{(\wnot^2-7\w^2)}=0
\end{equation}
Solving this equation for $\tan{\phi}$ and keeping only the leading term, the phase lag at amplitude resonance writes
\begin{equation}
    \tan{\phi_a} = \frac{3\w_a(\wnot^2-\w_a^2)}{\bar{\zeta}\wnot(\wnot^2-7\w_a^2)}
    \label{eq:phi_31}
\end{equation}
Assuming further that the ratio 
\begin{equation}
    \frac{\wnot^2-\w_a^2}{\wnot^2-7\w_a^2}\simeq 4
\end{equation}
yields:
\begin{equation}
    \tan{\phi_a} = \frac{12\w_a}{\bar{\zeta}\wnot}\label{PQ31}
\end{equation}
Inserting this relation in Equations \eqref{eq:drdphi31_0} and assuming that the static response is constant, i.e.,
$\Gamma\simeq\Gamma_*=\frac{9\gamma}{8\wnot}$, provides an expression of the amplitude of the third harmonic and of the frequency at amplitude resonance: 
\begin{equation}
    \begin{cases}
    A_{3,a} = \frac{\alpha\bar{\Gamma}^3_*}{2\bar{\zeta}\wnot^2\sqrt{\bar{\zeta}^2\wnot^2+144\w_a^2}}\\
    \w_a = \sqrt{\frac{-c_2+\sqrt{c_2^2-4c_1c_3}}{2c_1}}
    \end{cases}\label{ARC31a}
\end{equation}
where
\begin{equation}
    \begin{cases}
    c_1 = \frac{1728}{\alpha}\\
    c_2 = -144\left(2\bar{\Gamma}^{2}_*+ \frac{4\wnot^2}{3\alpha}-\frac{3\bar{\zeta}^2}{4\alpha}\right)\\
    c_3 = \left(\frac{2\bar{\zeta}^2\wnot^2}{3\alpha}-2\bar{\Gamma}^{2}_* - \frac{4\wnot^2}{3\alpha}\right)\bar{\zeta}^2\wnot^2 -  \frac{\alpha^2\bar{\Gamma}^{6}_*}{4\bar{\zeta}^2\wnot^2}
    \end{cases}\label{ARC31b}
\end{equation}

\subsubsection{Phase lag at phase resonance}

For weak to moderate damping, Equation (\ref{PQ31}) shows that amplitude resonance occurs near phase quadrature between the third harmonic of the displacement and the forcing. The phase resonance for the 3:1 superharmonic resonance can thus be associated with a phase lag of $\pi/2$. The averaged equations of motion \eqref{eq:drdphi31_0} become:
\begin{equation}
    \begin{cases}
    r_p = \frac{\alpha\Gamma^3}{24\zeta\wnot\w_p}\\
    r_p = \sqrt{\frac{4\W}{3\alpha}-2\Gamma^2}
    \end{cases}
    \label{eq:PRNM_31}
\end{equation}
If we assume again that $\Gamma\simeq\Gamma_*$, it is possible to derive a closed-form expression for $A_{3,p}$ and $\w_p$:
\begin{equation}
    \begin{cases}
    A_{3,p} = \frac{\alpha\bar{\Gamma}^3_*}{24\bar{\zeta}\wnot^2\w_p^2}\\
    \w_p = \sqrt{\frac{-c_2+\sqrt{c_2^2-4c_1c_3}}{2c_1}}
    \end{cases}\label{PRC31a}
\end{equation}
where
\begin{equation}
    \begin{cases}
    c_1 = \frac{1728}{\alpha}\\
    c_2 = -144\left(2\bar{\Gamma}^{2}_*+ \frac{4\wnot^2}{3\alpha}\right)\\
    c_3 = -\frac{\alpha^2\bar{\Gamma}^{6}_*}{4\bar{\zeta}^2\wnot^2}
    \end{cases}\label{PRC31b}
\end{equation}

\begin{figure}[htpb] 
  \begin{subfigure}[b]{0.5\linewidth}
    \centering
    \includegraphics[width=1\linewidth]{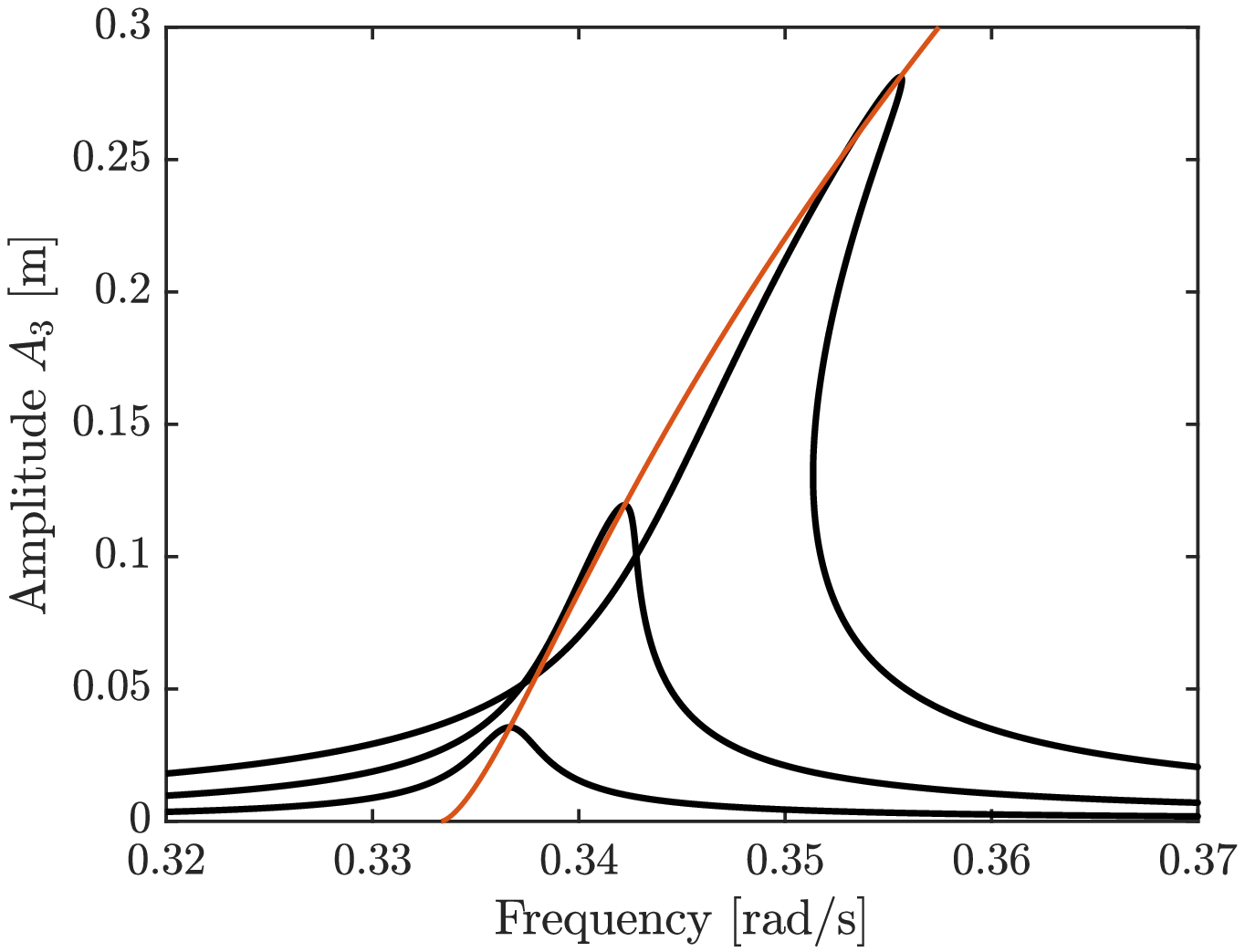} 
    \caption{\label{fig:31_NFRC}}
  \end{subfigure}
  \begin{subfigure}[b]{0.5\linewidth}
    \centering
    \includegraphics[width=1\linewidth]{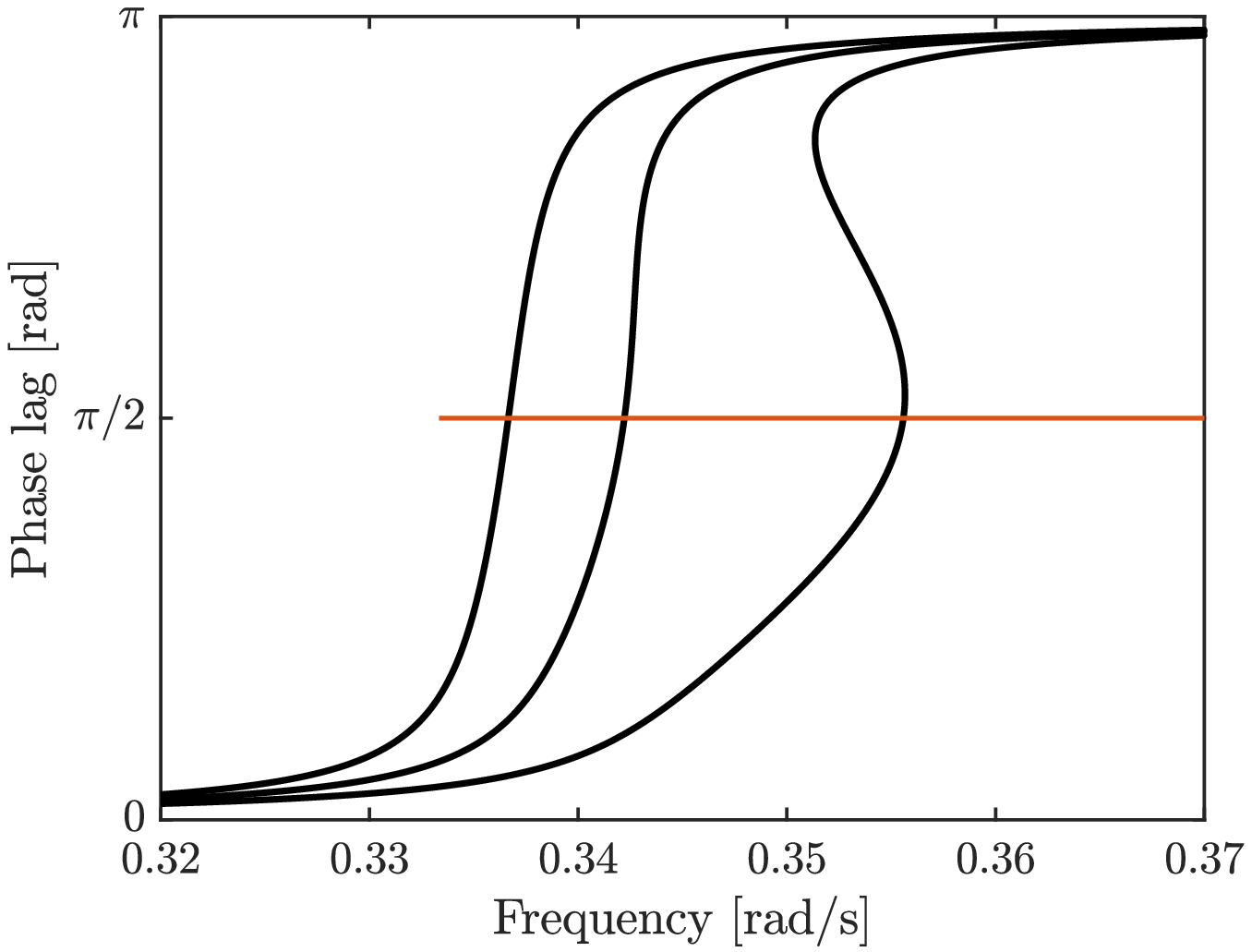}
    \caption{\label{fig:31_PHASE}}
  \end{subfigure} 
  \caption{Nonlinear frequency responses (black) and phase resonance/quadrature curves (orange) around the $3:1$ resonance of the Duffing oscillator for forcing amplitudes of 0.1N, 0.15N and 0.2N: (\subref{fig:31_NFRC}) amplitude and (\subref{fig:31_PHASE}) phase lag.}
  \label{fig:31RES_H3}
\end{figure}

Figures \ref{fig:31_NFRC} and \ref{fig:31_PHASE} compare the nonlinear frequency responses calculated from Equations \eqref{eq:drdphi31_0} and the phase resonance curves constructed thanks to Equations \eqref{PRC31a} and \eqref{PRC31b}. Clearly, the newly-defined concept of a phase resonance for the 3:1 superharmonic resonance is in excellent agreement with the maxima of the third harmonic of the response, at least for the amount of damping considered herein, i.e., 0.5\%. Assuming small $\bar{\zeta}$, this observation is also confirmed analytically by the direct comparison between Equations \eqref{ARC31a}-\eqref{ARC31b} and \eqref{PRC31a}-\eqref{PRC31b}. 

In the case of a softening Duffing oscillator ($\alpha<0$), the phase lag $\phi_p$ should be adjusted to $\frac{3\pi}{2}$ in order to have a positive amplitude $A_{3,p}$. This phase lag is still consistent with Equation \eqref{PQ31}.

\subsection{\texorpdfstring{$1:3$}{} resonance}
For the $1:3$ superharmonic resonance, \textit{i.e.}, $k=1$ and $\nu=3$, and $d=1$, first-order averaging gives:
\begin{equation}
    \begin{cases}
    \dot{r}=-\varepsilon\left(\zeta\wnot r - \frac{9\alpha\Gamma}{8\w}r^2\sin{3\phi}\right)\\
    \dot{\phi}=-\varepsilon\left(\frac{3\alpha}{8\w}\left(3r^2+6\Gamma^2-\frac{4\W}{\alpha}\right)-\frac{9\alpha\Gamma}{8\w }r\cos{3\phi}\right)
    \end{cases}
    \label{eq:drdphi13}
\end{equation}
For steady-state solutions,
\begin{equation}
    \begin{cases}
    \zeta\wnot =\frac{9\alpha\Gamma}{8\w}r\sin{3\phi}\\
    \frac{3\alpha}{8\w}\left(3r^2+6\Gamma^2-\frac{4\W}{\alpha}\right)=\frac{9\alpha\Gamma}{8\w }r\cos{3\phi}
    \end{cases}
    \label{eq:drdphi13_0}
\end{equation}

\subsubsection{Phase lag at amplitude resonance}
Amplitude resonance occurs when $\frac{\partial r}{\partial \w} = \frac{\partial r}{\partial \phi} = 0$:
\begin{equation}
    \begin{cases}
    \frac{\partial r}{\partial\phi} = \frac{8\zeta\wnot}{9\alpha\Gamma\sin{3\phi}}\left(\left(1-\frac{2\w^2}{\wnot^2-\w^2}\right)\frac{\partial \w}{\partial\phi}-\frac{3\w}{\tan{3\phi}}\right)=0\\
    
    \frac{\partial r}{\partial\w} = \frac{8\zeta\wnot}{9\alpha\Gamma\sin{3\phi}}\left(\left(1-\frac{2\w^2}{\wnot^2-\w^2}\right)-\frac{3\w}{\tan{3\phi}}\frac{\partial \phi}{\partial\w}\right)=0
    \end{cases}
    \label{eq:drdphi_drdw_13}
\end{equation}
We must have $\sin{3\phi}\neq 0$, \textit{i.e.}, $\phi\neq \frac{i\pi}{3}$, where $i$ is an integer. Following the same procedure as for the previous resonances gives:
\begin{equation}
    \frac{\partial r}{\partial\phi} = \frac{8\zeta\wnot}{9\alpha\Gamma\sin{3\phi}}\frac{\left(1-\frac{2\w^2}{\wnot^2-\w^2}\right)\varepsilon9\zeta\wnot - \frac{3\w}{\tan{3\phi}}\left(1-\varepsilon\frac{27\alpha\Gamma^2}{\wnot^2-\w^2}+\varepsilon\frac{6\zeta\wnot}{\wnot^2-\w^2}\frac{1}{\tan{3\phi}}\right)}{1-\varepsilon\frac{27\alpha\Gamma^2}{\wnot^2-\w^2}+\varepsilon\frac{6\zeta\wnot}{\wnot^2-\w^2}\frac{1}{\tan{3\phi}} - \left(1-\frac{2\w^2}{\wnot^2-\w^2}\right)\left(\varepsilon\frac{16\zeta^2\wnot^2}{3\alpha\Gamma^2\sin^2{3\phi}}+\varepsilon\frac{3\zeta\wnot}{\w\tan{3\phi}}\right)}
\end{equation}
The numerator is equal to 0 when:
\begin{equation}
    9\varepsilon  \zeta\wnot\left(1-\frac{2\w^2}{\wnot^2-\w^2}\right)\tan^2{3\phi}-3\w\left(1-\varepsilon\frac{27\alpha\Gamma^2}{\wnot^2-\w^2}\right)\tan{\phi} - \varepsilon\frac{18\zeta\wnot\w^2}{\wnot^2-\w^2}=0
\end{equation}
Solving this equation for $\tan{3\phi}$ and keeping only the leading term, the phase lag at amplitude resonance can be approximated with
\begin{equation}
    \tan{3\phi_a} = \frac{\w_a}{3\zeta\wnot\left(1 + \frac{2\w_a^2}{\w_a^2-\wnot^2}\right)}
    \label{eq:phi_13}
\end{equation}
Assuming further that the ratio 
\begin{equation}
    \frac{2\w_a^2}{\w_a^2-\wnot^2}\simeq \frac{9}{4}
\end{equation}
yields:
\begin{equation}
    \tan{3\phi_a} = \frac{4\w_a}{39\bar{\zeta}\wnot}
    \label{PQ13}
\end{equation}
Inserting this relation in Equations \eqref{eq:drdphi13_0} gives:
\begin{equation}
    \begin{cases}
    A_{1,a} = \frac{2\bar{\zeta}\wnot}{9\alpha\bar{\Gamma}}\sqrt{1521\bar{\zeta}^2\wnot^2+16\w_a^2}\\
    \bar{\gamma} = \frac{\|\wnot^2-\w_a^2\|}{\sqrt{6\alpha}}\sqrt{(2\bar{\W}+13\bar{\zeta}^2\wnot^2)^2\pm\sqrt{(2\bar{\W}+13\bar{\zeta}^2\wnot^2)^2-\frac{8}{9}\bar{\zeta}^2\wnot^2(1521\bar{\zeta}^2\wnot^2+16\w_a^2)}}
    \end{cases}\label{RA13}
\end{equation}
where $\bar{\W}=\w_k^2-\wnot^2$. Unlike the $3:1$ superharmonic resonance, the static response cannot be assumed to be constant because the frequency varies much faster for the 1:3 subharmonic resonance (see Figure \ref{fig:13_NFRC}). An explicit expression for the resonance frequency $\w_a$ as a function of the forcing $\bar{\gamma}$ can thus not be derived. We also note that, due to the $\pm$ sign, there exist two frequencies satisfying \eqref{RA13}, the greatest (lowest) frequency corresponding to the maximum (minimum) response on the isolated branch. It is thus the greatest frequency which is in relation with the resonance frequency $\w_a$.  

\subsubsection{Phase lag at phase resonance}

For weak to moderate damping, Equation (\ref{PQ13}) shows that amplitude resonance occurs near phase lags equal to  $\frac{\pi}{6}+\frac{i\pi}{3}$ where $i$ is an integer. For odd (even) values of $i$, $r$ takes positive (negative) values. Considering positive amplitudes, the phase resonance for the 1:3 subharmonic resonance can be associated with phase lags equal to $\frac{\pi}{2}, \frac{7\pi}{6}\,\mbox{and}\,\frac{11\pi}{6}$.  For $\frac{\pi}{2}$, the averaged equations of motion \eqref{eq:drdphi31_0} can be transformed into:
\begin{equation}
    \begin{cases}
    A_{1,p} = \frac{8\bar{\zeta}\wnot\w_p}{9\alpha\bar{\Gamma}}\\
    \bar{\gamma} = \frac{\|\wnot^2-\w_p^2\|}{\sqrt{3\alpha}}\sqrt{\bar{\W}\pm\sqrt{\bar{\W}^2-\frac{32}{9}\bar{\zeta}^2\wnot^2\w_p^2}}\label{PRC13b}
    \end{cases}
\end{equation}
The same expressions can be obtained if the two other phase lags are considered instead.

\begin{figure}[htpb] 
  \begin{subfigure}[b]{0.5\linewidth}
    \centering
    \includegraphics[width=1\linewidth]{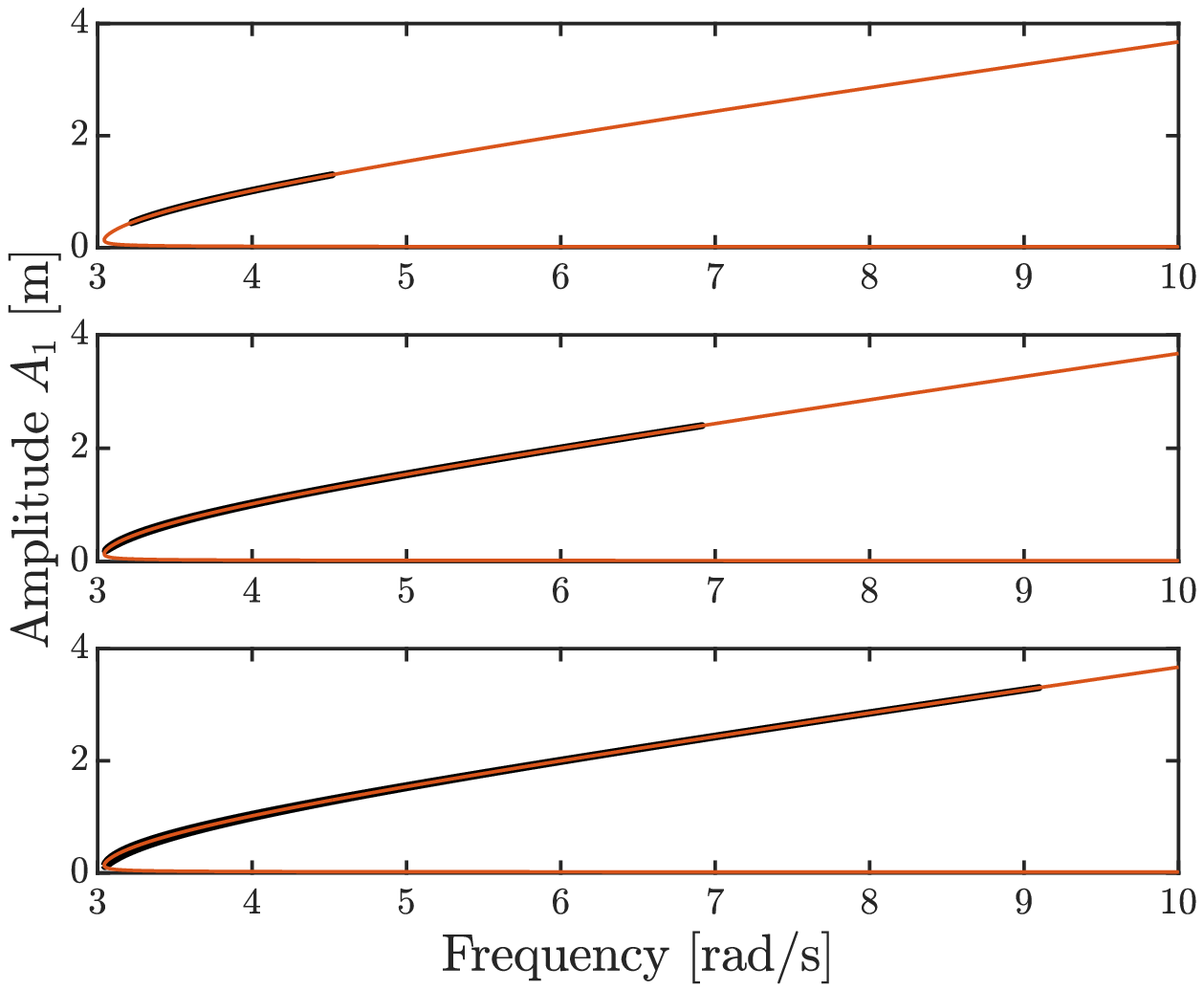} 
    \caption{\label{fig:13_NFRC}}
  \end{subfigure}
  \begin{subfigure}[b]{0.5\linewidth}
    \centering
    \includegraphics[width=1\linewidth]{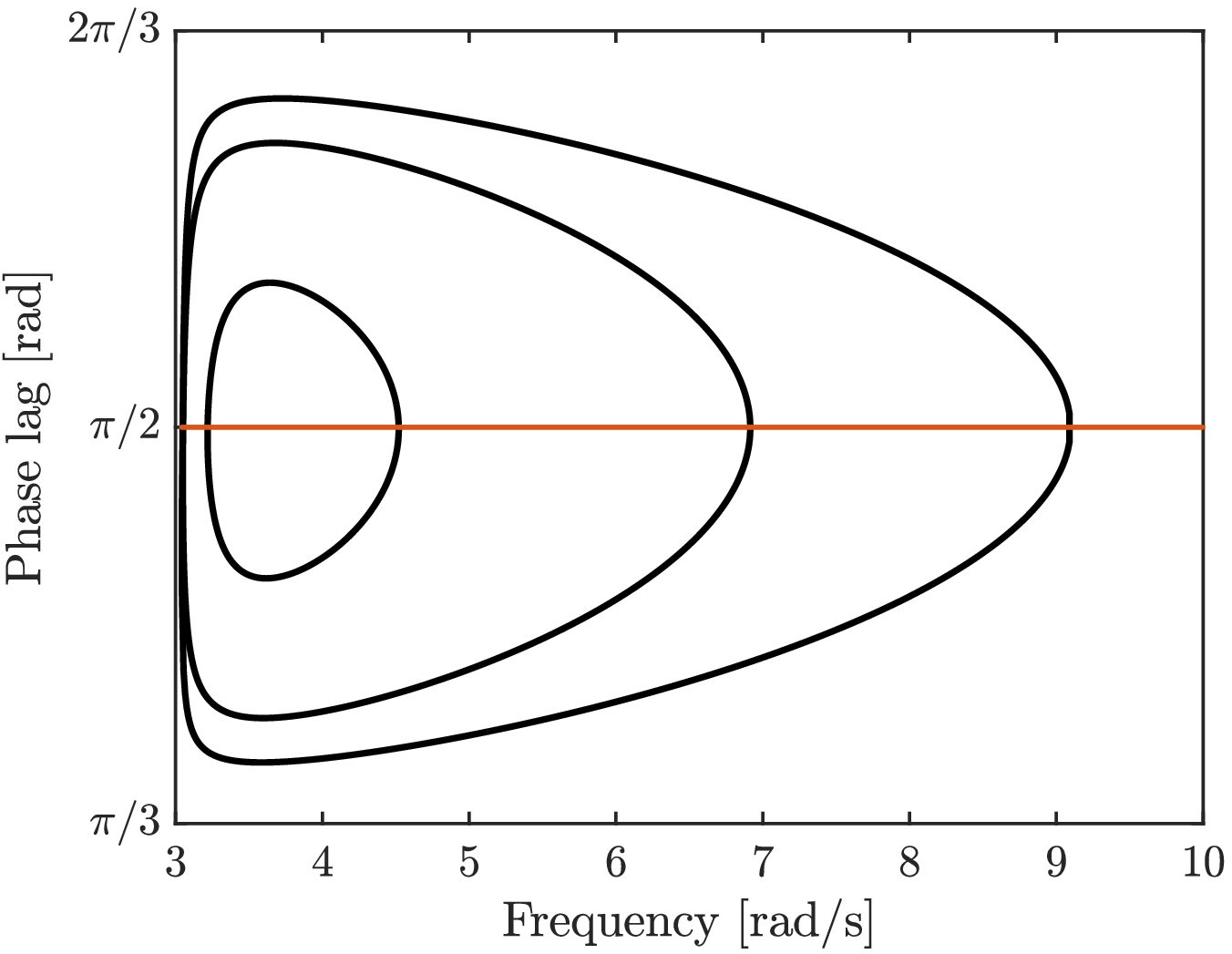}
    \caption{\label{fig:13_PHASE}}
  \end{subfigure} 
  \caption{Nonlinear frequency responses (black) and phase resonance/quadrature curves (orange) around the $1:3$ resonance of the Duffing oscillator for forcing amplitudes of 0.3N, 0.6N and 1N: (\subref{fig:13_NFRC}) amplitude and (\subref{fig:13_PHASE}) phase lag.}
  \label{fig:13RES_H1}
\end{figure}

Figures \ref{fig:13_NFRC} and \ref{fig:13_PHASE} compare the nonlinear frequency responses calculated from Equations \eqref{eq:drdphi13_0} and the phase resonance curve constructed numerically thanks to Equations \eqref{PRC13b}. The phase quadrature curve is found to trace out the locus of the maxima of the different isolated responses.

For a softening Duffing oscillator, amplitude resonance still occurs for phase lags $\phi_p$ near $\frac{\pi}{6}+\frac{i\pi}{3}$ except that positive amplitudes occur now when $i$ is odd. Thus, the resonant phase lags are $\frac{\pi}{6}, \frac{5\pi}{6}\,\mbox{and}\,\frac{3\pi}{2}$.

\subsection{\texorpdfstring{$1:2$}{} resonance}
Using Equation \eqref{eq:EOM1D_SCALE} with $d=2$, second-order averaging yields:
\begin{equation}
    \begin{cases}
    \dot{r} = -\frac{\varepsilon^2}{2}\left(2\zeta\wnot r + \frac{33\alpha^2\Gamma^2}{4\w^3}r^3\sin{4\phi}\right)\\
    \dot{\phi} = -\frac{\varepsilon\alpha}{4\w}\left(3r^2+6\Gamma^2-\frac{4\W}{\alpha}\right) + \frac{\varepsilon^2}{2}\left(R_{1:2}(r^2)-\frac{33\alpha^2\Gamma^2}{4\w^3}r^2\cos{4\phi}\right)
    \end{cases}
\end{equation}
where 
\begin{equation}
    R_{1:2}(r^2)=\left(\frac{2\W^2}{\w^3}-\frac{6\alpha\Gamma^2\W}{\w^3}-\frac{51\alpha^2\Gamma^2}{10\w^3}\right)r^2 - \left(\frac{6\alpha\W}{\w^3}+\frac{33\alpha^2\Gamma^2}{4\w^3}\right)r^4+\frac{51\alpha^2}{16\w^3}r^6
\end{equation}
Steady-state solutions obey:
\begin{equation}
    \begin{cases}
    2\zeta\wnot =  - \frac{33\alpha^2\Gamma^2}{4\w^3}r^2\sin{4\phi}\\
    \frac{\alpha}{2\w}\left(3r^2+6\Gamma^2-\frac{4\W}{\alpha}\right) =  \varepsilon\left(R_{1:2}(r^2)-\frac{33\alpha^2\Gamma^2}{4\w^3}r^2\cos{4\phi}\right)
    \end{cases}
    \label{eq:drdphi12_0}
\end{equation}
\subsubsection{Phase lag at amplitude resonance}
Neglecting the $\mathcal{O}(\varepsilon)$ term in the second Equation of \eqref{eq:drdphi12_0} gives an approximation $r_0$ of the amplitude $r$:
\begin{equation}
    r_0 = \sqrt{\frac{4\W}{3\alpha}-2\Gamma^2}
    \label{eq:r0}
\end{equation}
Its derivative is:
\begin{equation}
    \frac{\partial r_0}{\partial \w} = \frac{4}{r_0}\left(\frac{1}{12\varepsilon\alpha}-\frac{\gamma^2}{(\wnot^2-\w^2)^3}\right)\w\label{der}
\end{equation}

\begin{figure}[htpb] 
    \centering
    \includegraphics[width=1\linewidth]{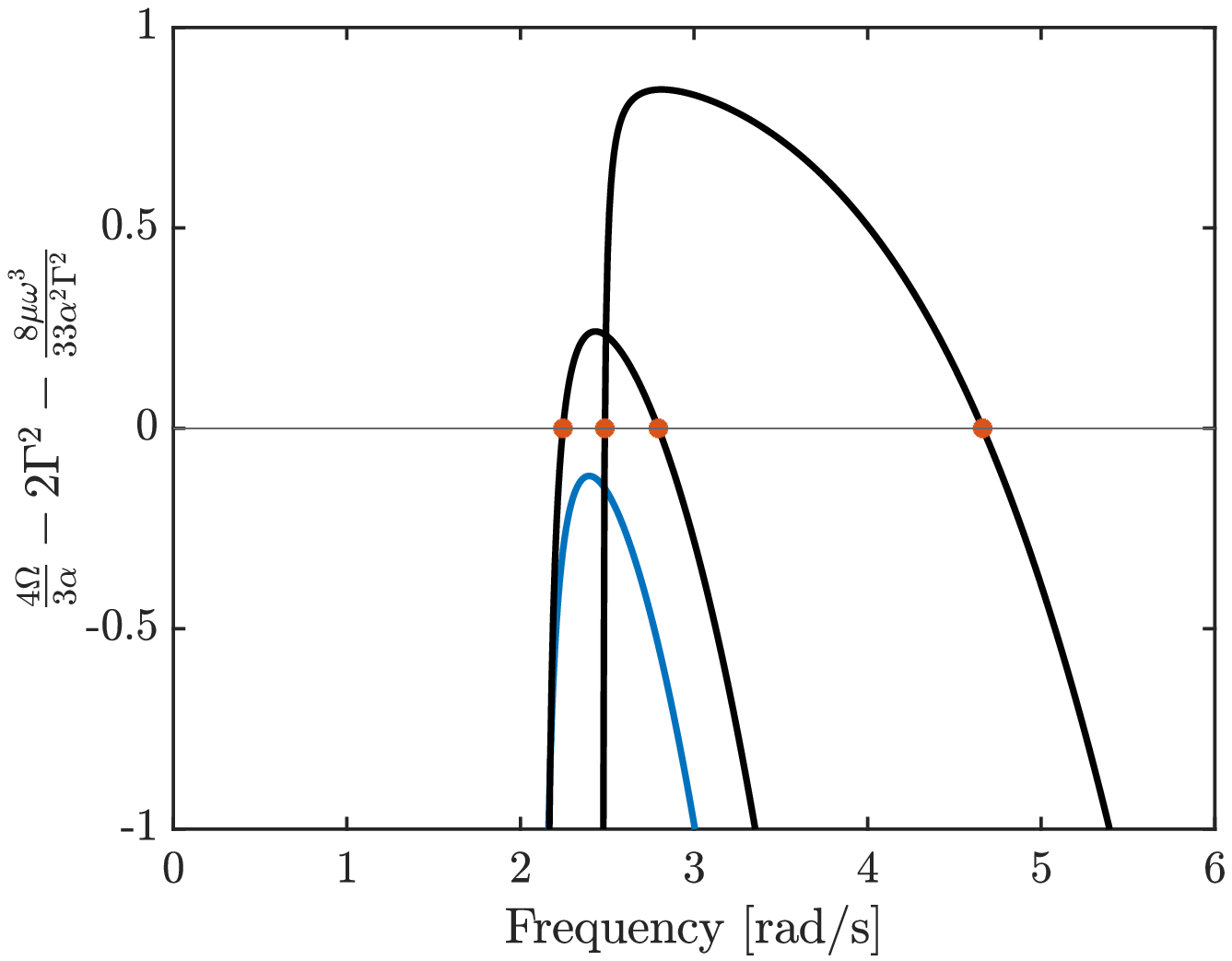}
    \caption{\label{fig:12_ineq_w}}
  \caption{Numerical verification of the inequality  \eqref{eq:r0_existence2_12} for a forcing of 0.8N (blue), 1N and 3N (black).}
  \label{fig:12_rangew}
\end{figure}

Considering that the sinus function in the first equation of \eqref{eq:drdphi12_0} is bounded by $-1$ and $1$, an existence condition for $r$ is derived:
\begin{equation}
    -1\leq -\frac{8\zeta\wnot\w^3}{33\alpha^2\Gamma^2r^2} \leq 1
\end{equation}
The second inequality is always true. $r_0$ is thus injected in the first inequality:
\begin{equation}
    \frac{4\W}{3\alpha}\geq 2\Gamma^2+\frac{8\zeta\wnot\w^3}{33\alpha^2\Gamma^2}
    \label{eq:r0_existence2_12}
\end{equation}
The numerical resolution in Figure \ref{fig:12_ineq_w} indicates that, if the forcing exceeds a certain threshold, there exist two frequencies, $\w_{inf}$ and $\w_{sup}$, which define the domain of existence of the 1:2 subharmonic resonance. Conversely, if the forcing is too low, the inequality is not satisfied, and the 1:2 subharmonic resonance does not exist. Because Equation \eqref{der} shows that $r_0$ is increasing monotonically with respect to frequency since $\alpha>0$, $r_0$ is thus maximum (minimum) when $\w$ is equal to $\w_{sup}$ ($\w_{inf}$), and amplitude resonance occurs when $\omega=\w_{sup}$. 

\subsubsection{Phase lag at phase resonance}

For this specific resonance, because the $\mathcal{O}(\varepsilon)$ term was neglected in Equations \eqref{eq:drdphi12_0}, phase resonance corresponds to amplitude resonance and takes place when $\omega=\w_{sup}$ or, equivalently, $\sin{4\phi_a}=-1$. Phase resonance is thus defined for a phase lag $\phi_a=\frac{3\pi}{8}+\frac{i\pi}{2}$, where $i=0,1,2,3$, which transforms \eqref{eq:drdphi12_0} into:
\begin{equation}
    \begin{cases}
    r_p = \sqrt{\frac{8\zeta\wnot\w^3}{33\alpha^2\Gamma^2}}\\
    r_p = \sqrt{\frac{4\W}{3\alpha}-2\Gamma^2}
    \end{cases}
    \label{eq:PRNM_12}
\end{equation}
from which the forcing $\gamma$ can be computed as a function of the resonant frequency $\w_p$:
\begin{equation}
    \begin{cases}
    A_p = \sqrt{\frac{8\bar{\zeta}\wnot\w_p^3}{33\alpha^2\bar{\Gamma}^2}}\\
    \bar{\gamma} = \frac{\|\wnot^2-\w_p^2\|}{\sqrt{3\alpha}}\sqrt{\bar{\W}\pm\sqrt{\bar{\W}^2-\frac{12}{11}\bar{\zeta}\wnot\w_p^3}}
    \end{cases}\label{PRC12}
\end{equation}
The $\pm$ sign indicates that, for a fixed forcing, there exist two possible frequencies corresponding to the minimum and maximum values of $A_p$. 

\begin{figure}[htpb] 
  \begin{subfigure}[b]{0.5\linewidth}
    \centering
    \includegraphics[width=1\linewidth]{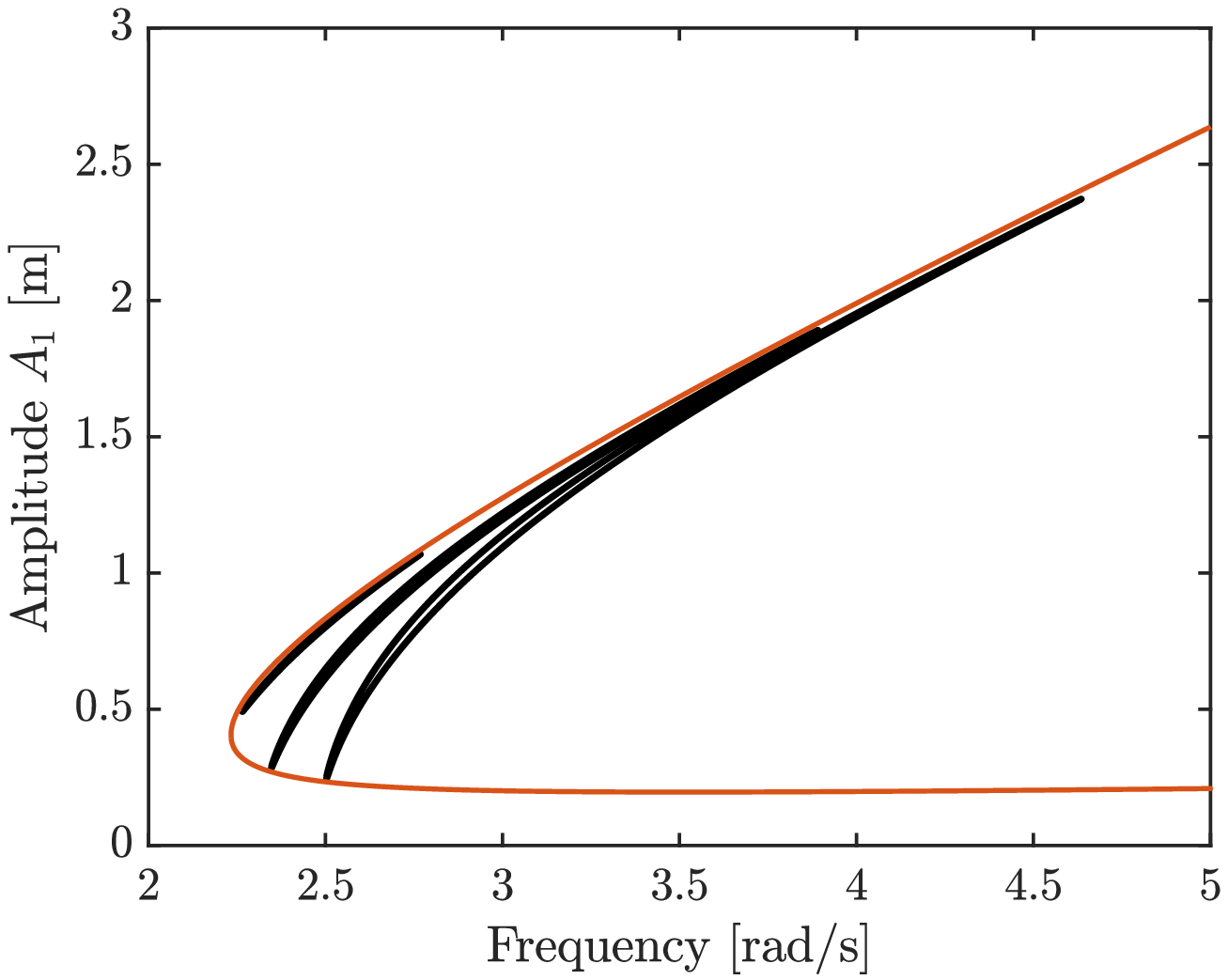} 
    \caption{\label{fig:12_NFRC}}
  \end{subfigure}
  \begin{subfigure}[b]{0.5\linewidth}
    \centering
    \includegraphics[width=1\linewidth]{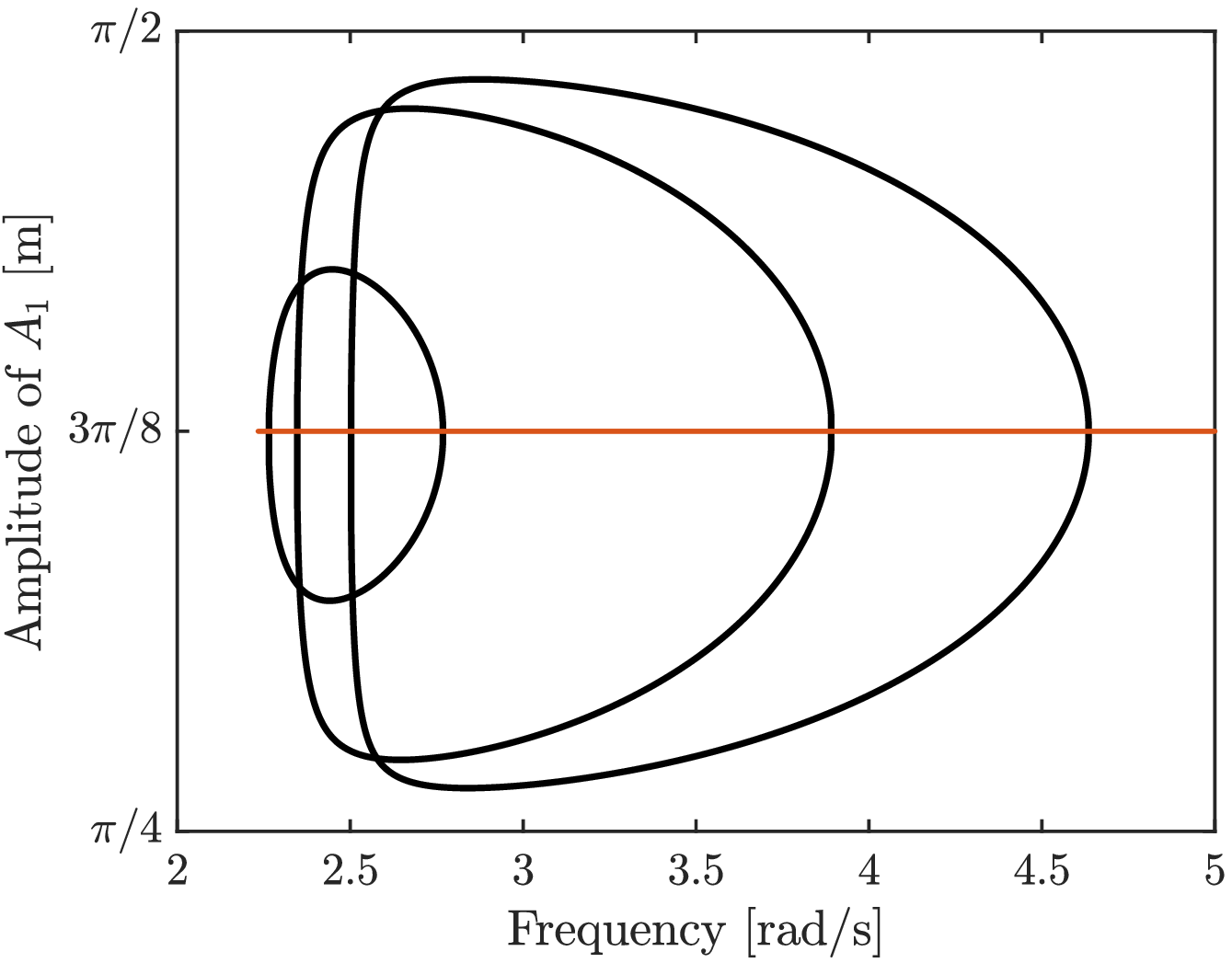}
    \caption{\label{fig:12_PHASE}}
  \end{subfigure} 
  \caption{Nonlinear frequency responses (black) and phase resonance curves (orange) around the $1:2$ resonance of the Duffing oscillator for forcing amplitudes of 1N, 2N and 3N: (\subref{fig:12_NFRC}) amplitude and (\subref{fig:12_PHASE}) phase lag.}
  \label{fig:12RES_H1}
\end{figure}

Figures \ref{fig:12_NFRC} and \ref{fig:12_PHASE} compare the nonlinear frequency responses calculated from Equations \eqref{eq:drdphi12_0} and the phase resonance curve constructed numerically thanks to Equations \eqref{PRC12}. The phase resonance curve is found to trace out the locus of the maxima and minima of the different isolated branches.

    \section{Additional Higher-order Resonances}
Higher-order averaging is required to obtain the resonant motions around other $k:\nu$ resonances. For each resonance, the second relation of the governing equations takes the form:
\begin{equation}
    3r^2+6\Gamma^2-\frac{4\W}{\alpha} = \mathcal{O}(\varepsilon)
\end{equation}
Neglecting the $\mathcal{O}(\varepsilon)$ term, an approximation $r_0$ of the amplitude $r$ is obtained as in Equation \eqref{eq:r0}. The corresponding derivative is:
\begin{equation}
    \frac{\partial r_0}{\partial \w} = \frac{4}{r_0}\left(\frac{k^2}{3\nu^2\varepsilon\alpha}-\frac{\gamma^2}{(\wnot^2-\w^2)^3}\right)\w
    \label{eq:dr0dw}
\end{equation}
which is always positive when $\w>\wnot$, i.e., when $\nu>k$. If $\w<\wnot$, the derivative can be either positive, negative or positive then negative depending on the value of $\w$ and $\gamma$.

\subsection{Subharmonic resonances (\texorpdfstring{$1:\nu$}{})}
\subsubsection{Odd subharmonic resonances}
These resonances occur when $\nu$ takes odd values. Considering first the $1:5$ resonance, the governing equations at steady state are:
\begin{equation}
    \begin{cases}
    2\zeta\wnot =  - \frac{1875\alpha^2\Gamma}{128\w^3}r^3\sin{5\phi}\\
    3r^2+6\Gamma^2-\frac{4\W}{\alpha} = \mathcal{O}(\varepsilon)
    \end{cases}
    \label{eq:drdphi15_0}
\end{equation}
In this case, $\frac{\partial r_0}{\partial \w}$ is always positive. Applying the same reasoning as for the $1:2$ resonance indicates that phase resonance occurs for a phase lag equal to $\frac{\pi}{10}+\frac{2i\pi}{5}$, where $i=0,1,2,3,4$. Because all these phase lags provide the same results, $\frac{\pi}{2}$ can thus be chosen as the resonant phase lag. This finding can be extended to higher odd values of $\nu$.

\subsection{Even subharmonic resonances}
The governing equations of the $1:4$ resonance are:
\begin{equation}
    \begin{cases}
    \zeta\wnot = -\frac{1665\alpha^4\Gamma^2}{\w^6}r^6\sin{8\phi}\\
    3r^2+6\Gamma^2-\frac{4\W}{\alpha} = \mathcal{O}(\varepsilon)
    \end{cases}
\end{equation}
$\frac{\partial r_0}{\partial \w}$ is again positive. Phase resonance occurs for a phase lag equal to $\frac{3\pi}{16}+\frac{2i\pi}{4}$, where $i=0,1,2,3$, meaning that $\frac{3\pi}{16}$ can be chosen as the resonant phase lag. For higher even values of $\nu$, the value $\frac{3\pi}{4\nu}$ is to be considered.

\subsection{Superharmonic resonance (\texorpdfstring{$k:1$}{})}
\subsubsection{Odd superharmonic resonances}
For the $5:1$ resonance, one has:
\begin{equation}
    \begin{cases}
    2\zeta\wnot r =  \frac{3\alpha^2\Gamma^5}{1280\w^3}\sin{\phi}\\
    3r^2+6\Gamma^2-\frac{4\W}{\alpha} =  \mathcal{O}(\varepsilon)
    \end{cases}
    \label{eq:drdphi51_0}
\end{equation}
Because the sine function is bounded by $-1$ and $1$ and because the argument of the square root inside $r_0$ must be positive, the inequality
\begin{equation}
    2\Gamma^2\leq\frac{4\W}{3\alpha} \leq \left(\frac{3\alpha^2\Gamma^5}{2560\zeta\wnot\w^3}\right)^2+2\Gamma^2
\end{equation}
must be verified. In the domain of existence, it may happen that $\frac{\partial r_0}{\partial \w}$ be negative, but this only occurs for large values of $\gamma$ or when $\w$ is close to $\wnot$, which is not consistent with the hypotheses of the averaging technique. We can thus consider that $r_0$ increases monotonically and that $\frac{\pi}{2}$ is the phase lag which characterizes both the amplitude and phase resonances. This also holds for higher-order odd superharmonic resonances.

\subsubsection{Even superharmonic resonances}
Applying second-order averaging around the $2:1$ resonance yields:
\begin{equation}
    \begin{cases}
    2\zeta\wnot =  - \frac{21\alpha^2\Gamma^4}{640\w^3}\sin{2\phi}\\
    3r^2+6\Gamma^2-\frac{4\W}{\alpha} = \mathcal{O}(\varepsilon)
    \end{cases}
    \label{eq:drdphi21_0}
\end{equation}
Contrary to other resonances, the first relation of \eqref{eq:drdphi21_0} shows that there is no direct relation between $r$ and $\phi$. In fact, numerical and analytical solutions for the even superharmonic resonances do not match. We thus conclude that our developments cannot be  accurately predict these resonances.

\subsection{Ultra-subharmonic resonances (\texorpdfstring{$k:\nu$}{})}
\subsubsection{\texorpdfstring{$\nu>k$}{}}
Fourth-order averaging is used to derive the governing equations for the $2:3$ ultra-subharmonic resonance.
\begin{equation}
    \begin{cases}
    24\zeta\wnot =  - \frac{297257881995}{282591232}\frac{\alpha^4\Gamma^4}{\w^7}r^4\sin{6\phi}\\
    3r^2+6\Gamma^2-\frac{4\W}{\alpha} = \mathcal{O}(\varepsilon)
    \end{cases}
\end{equation}
Amplitude and phase resonances occur for a phase lag equal to $\frac{\pi}{4}+\frac{i\pi}{3}$ where $i=0,1,\dotso,5$. 

\subsubsection{\texorpdfstring{$k>\nu$}{}}
\begin{figure}[htpb] 
  \begin{subfigure}[b]{0.5\linewidth}
    \centering
    \includegraphics[width=1\linewidth]{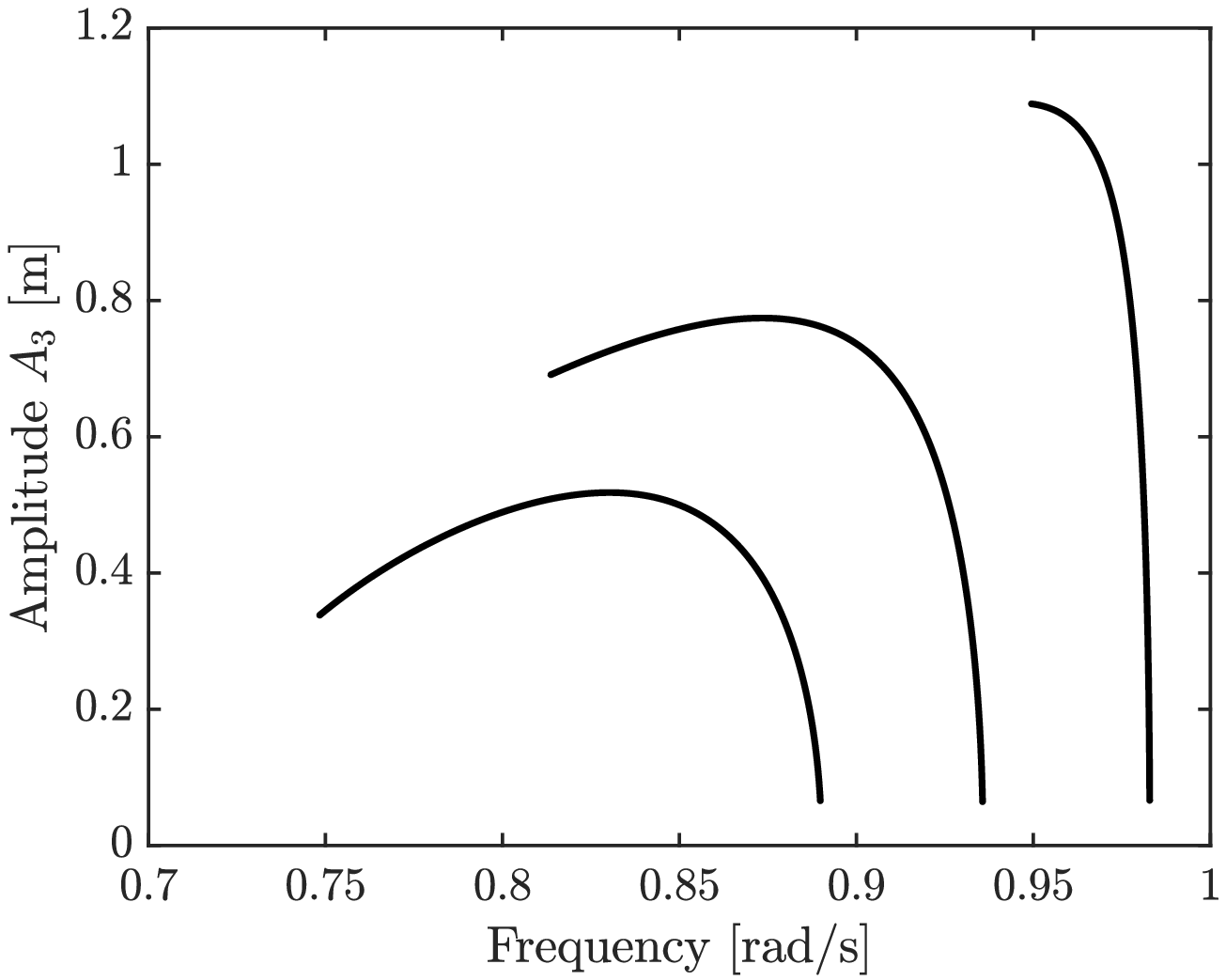} 
    \caption{\label{fig:32_NFRC}}
  \end{subfigure}
  \begin{subfigure}[b]{0.5\linewidth}
    \centering
    \includegraphics[width=1\linewidth]{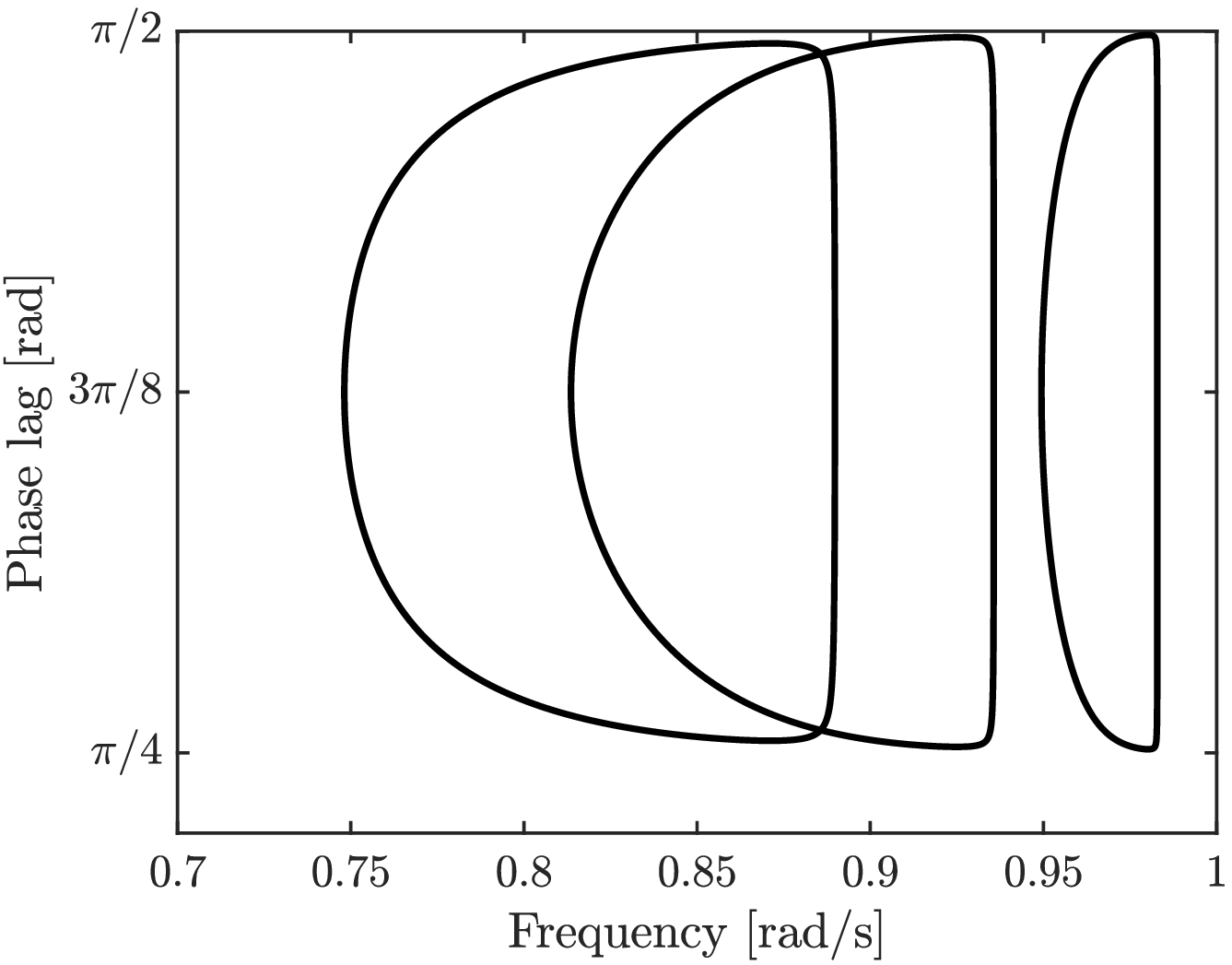}
    \caption{\label{fig:32_PHASE}}
  \end{subfigure} 
  \caption{Nonlinear frequency responses around the $3:2$ resonance of the Duffing oscillator for forcing amplitudes of 0.03N, 0.1N and 0.15N (right to left in the graphs): (\subref{fig:32_NFRC}) amplitude and (\subref{fig:32_PHASE}) phase lag.}
  \label{fig:32RES_H3}
\end{figure}

For the $3:2$ resonance,
\begin{equation}
    \begin{cases}
    24\zeta\wnot =  - \frac{1973735}{22289904}\frac{\alpha^4\Gamma^6}{\w^7}r^2\sin{4\phi}\\
    3r^2+6\Gamma^2-\frac{4\W}{\alpha} = \mathcal{O}(\varepsilon)
    \end{cases}
\end{equation}
In this case, $r_0$ does not increase monotonically since  $\frac{\partial r_0}{\partial\w}$ may be negative. However, the numerical simulations in Figure \ref{fig:32RES_H3} show that the resonance branch exists between two frequencies which feature a phase lag equal to $\frac{3\pi}{8}+\frac{i\pi}{2}$. Even though this phase lag does no longer correspond to amplitude resonance, it still provides valuable information, because it locates the extremities of the isolated response.

    \section{Numerical Validation}




The analytical results in the previous sections considered the amplitude and phase lag of the harmonic $k$, i.e., the harmonic triggering the $k:\nu$ resonance. We observed that amplitude resonance was occurring near a well-defined phase lag, allowing us to extend the concept of a phase resonance to secondary resonances. The phase resonances of the Duffing oscillator can be classified into two families depending on the value $\phi_k$:
\begin{itemize}
    \item $\phi_k=\frac{\pi}{2}$ (phase quadrature) when $k$ and $\nu$ are odd;
    \item $\phi_k=\frac{3\pi}{4\nu}$ when either $k$ or $\nu$ is even.
\end{itemize}
Though the averaging technique did not give satisfying results for the even superharmonic resonances, their resonant phase lag follows the above rule, i.e., they take the value $\frac{3\pi}{4}$, as evidenced numerically in \cite{JSVVOLVERT}.

\begin{figure}[htpb] 
  \begin{subfigure}[b]{0.5\linewidth}
    \centering
    \includegraphics[width=1\linewidth]{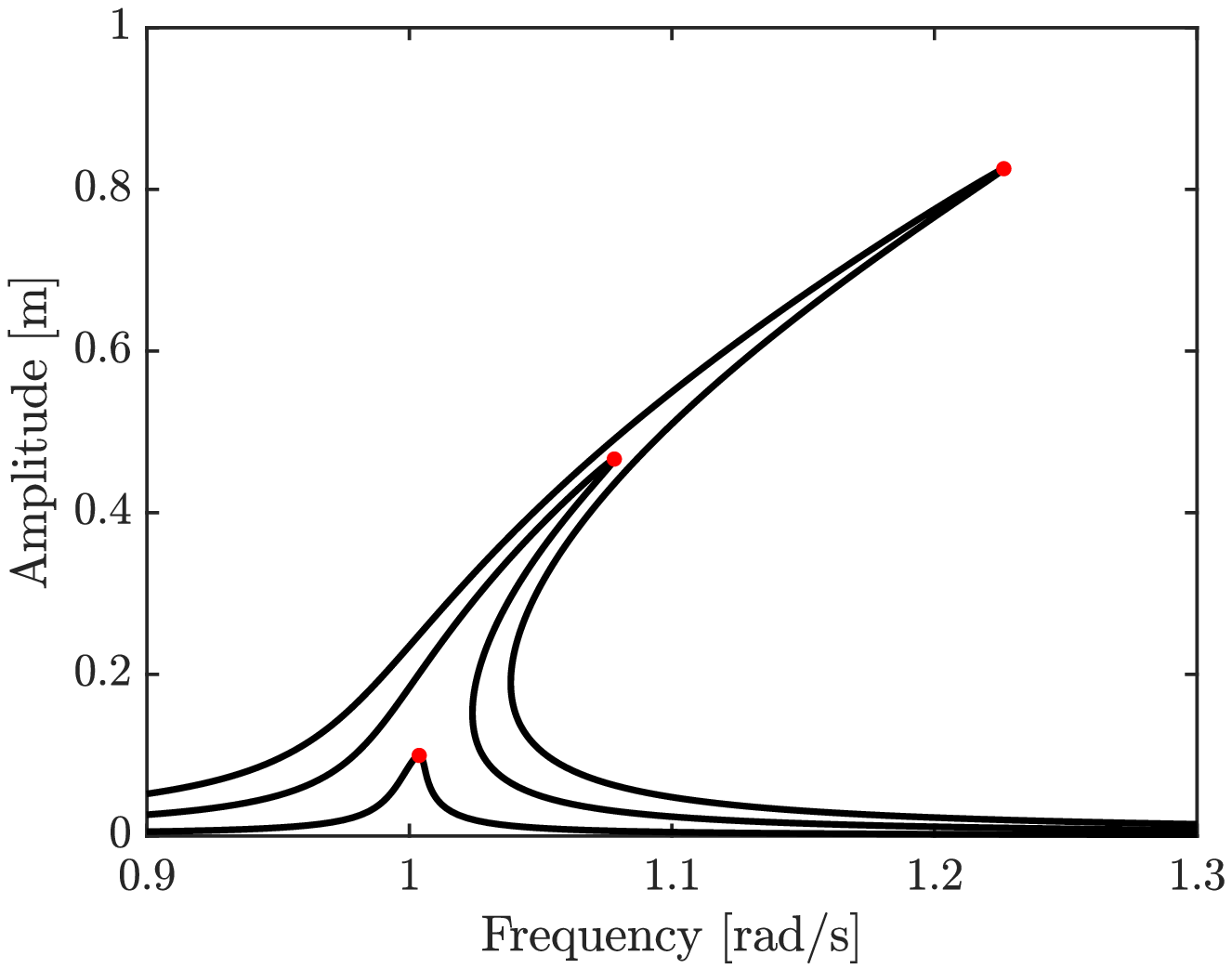} 
    \caption{\label{fig:FUNDRED_NFRC_NO_PRNM}}
  \end{subfigure}
  \begin{subfigure}[b]{0.5\linewidth}
    \centering
    \includegraphics[width=1\linewidth]{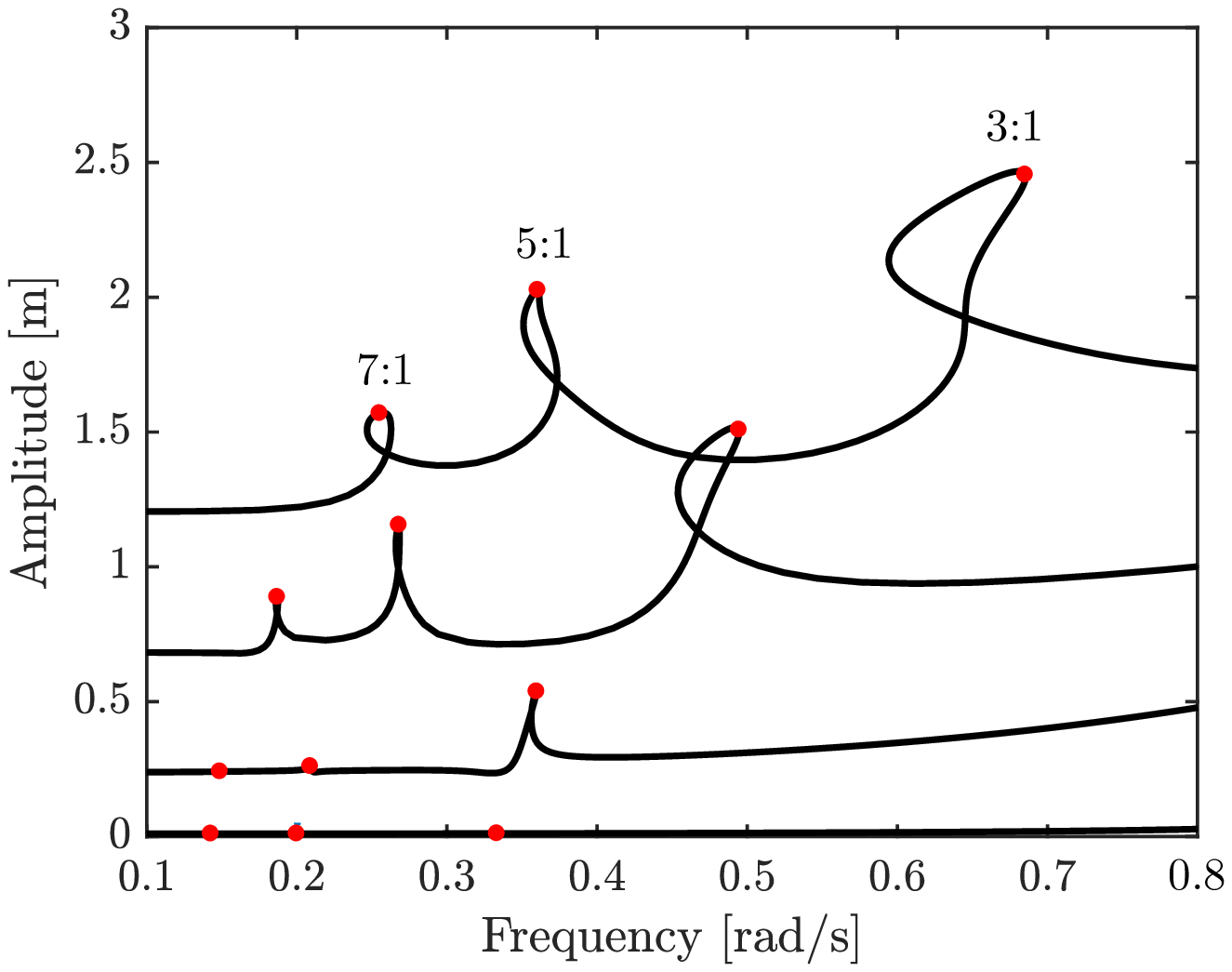} 
    \caption{\label{fig:SUPERHODD_NFRC_NO_PRNM}}
  \end{subfigure}
  
  \begin{subfigure}[b]{0.5\linewidth}
    \centering
    \includegraphics[width=1\linewidth]{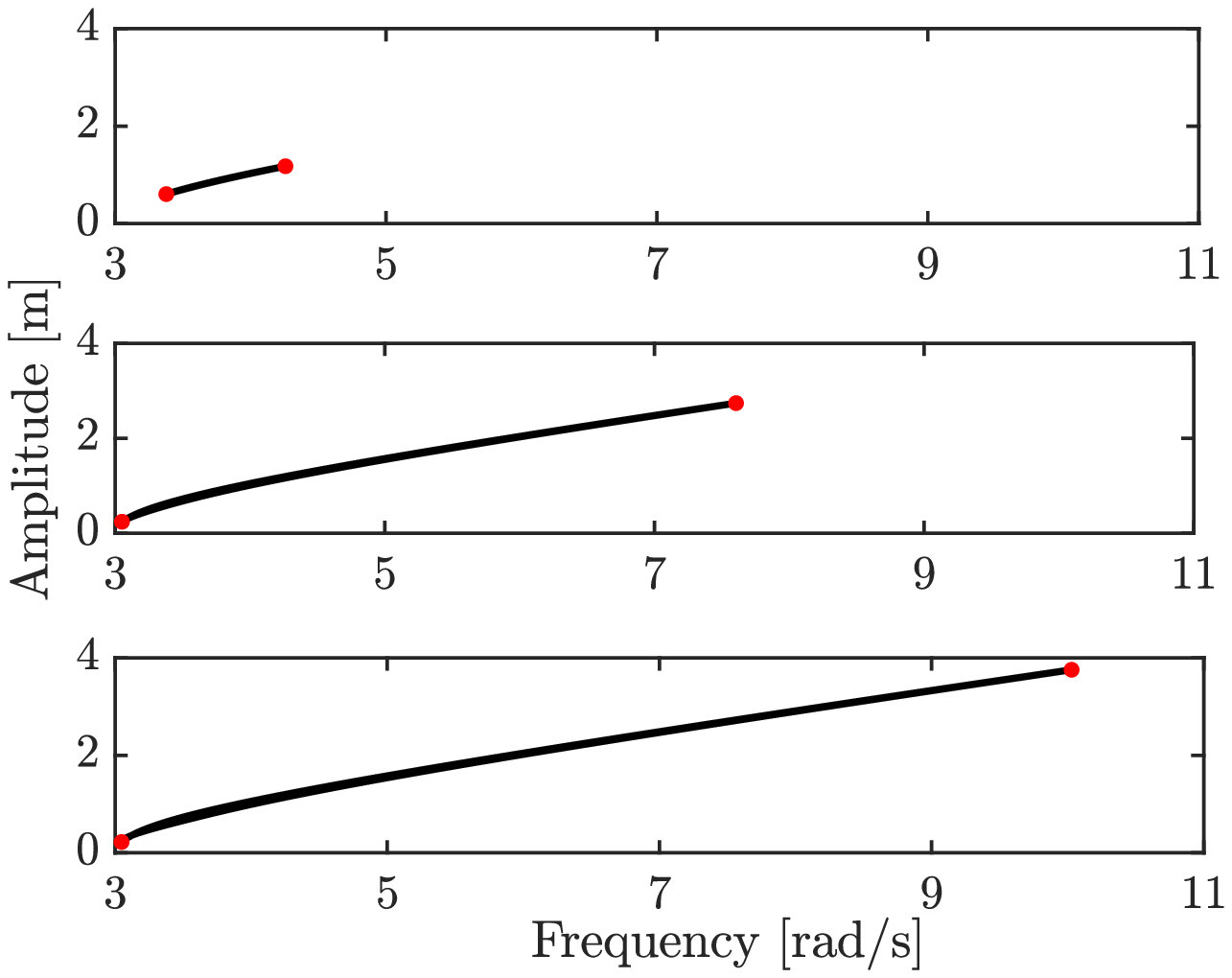}
    \caption{\label{fig:SUBHODD_NFRC_NO_PRNM}}
  \end{subfigure}
  \begin{subfigure}[b]{0.5\linewidth}
    \centering
    \includegraphics[width=1\linewidth]{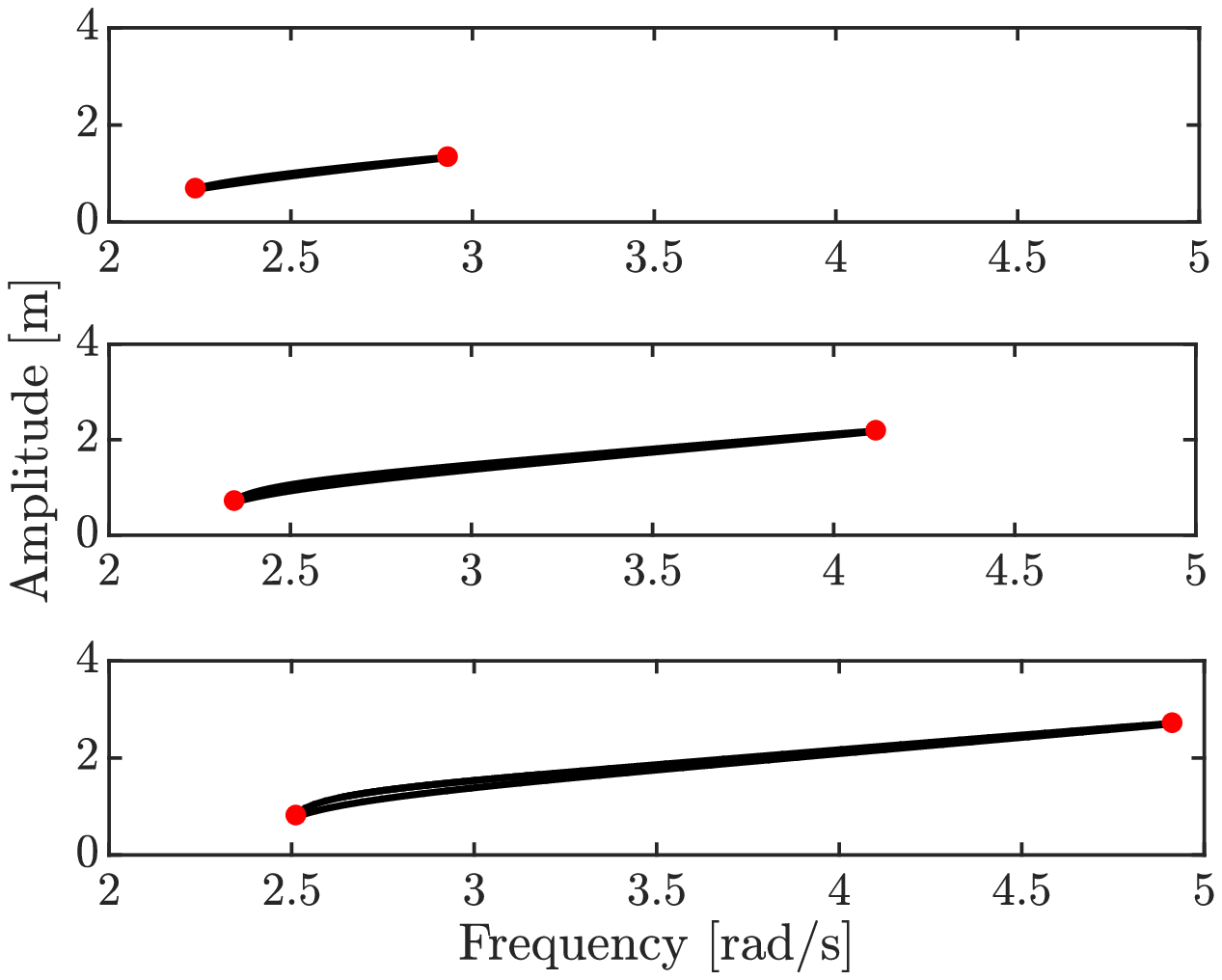}
    \caption{\label{fig:SUBHEVEN_NFRC_NO_PRNM}}
  \end{subfigure}
  \caption{NFRCs (black) and the corresponding phase resonance points (red) of the (\subref{fig:FUNDRED_NFRC_NO_PRNM}) $1:1$, (\subref{fig:SUPERHODD_NFRC_NO_PRNM}) $k:1$, (\subref{fig:SUBHODD_NFRC_NO_PRNM}) $1:3$ and (\subref{fig:SUBHEVEN_NFRC_NO_PRNM}) $1:2$ resonances.}
  \label{fig:PRNM_points} 
\end{figure}

To further validate the relevance of our developments, Figure \ref{fig:PRNM_points} represents the nonlinear frequency responses of the Duffing oscillator calculated numerically using the harmonic balance method \cite{DETROUX} for the $1:1$, $3:1$, $5:1$, $7:1$, $1:2$ and $1:3$ resonances. Unlike the previous figures, Figure \ref{fig:PRNM_points} depicts the multi-harmonic response of the Duffing oscillator. The red dots are located where phase resonance of the $k$-th harmonic occurs, \textit{i.e.}, $\phi_k=\frac{\pi}{2}$ for the $1:1$, $k:1$ and $1:3$ resonances and $\frac{3\pi}{8}$ for the $1:2$ resonance. We can clearly see that the so-defined phase resonance points can also accurately capture the amplitude resonance of the multi-harmonic response, and not only of the $k$-th harmonic.  

    \section{Conclusion}

The key contribution of this paper is the analytical characterization of the resonant phase lags of a hardening Duffing oscillator. For the $k:\nu$ resonance, the phase lag is computed between the $k$-th harmonic of the displacement and the harmonic forcing. When $k$ and $\nu$ are odd, phase resonance occurs when phase quadrature is achieved. When either $k$ or $\nu$ is even, phase resonance takes place for a phase lag equal to $3\pi/4\nu$. In almost all cases considered, phase resonance appears in the immediate vicinity of the amplitude resonance of the $k$th harmonic, at least for the amount of damping considered in this study.

These analytical results are in complete agreement with the numerical observations made in \cite{JSVVOLVERT}. They thus confirm the relevance of the concept of a phase resonance nonlinear mode (PRNM) which was defined as the point on the nonlinear frequency response which fulfills the phase resonance conditions. Eventually, the two papers lay down the foundations for rigorous phase resonance testing of nonlinear systems using phase-locked loops \cite{PETER} and the subsequent correlation between numerical and experimental analyses.

Future work should generalize these results to other types on nonlinearities, including softening nonlinearities, and to higher-dimensional systems.

    \bibliographystyle{unsrt} 
    \bibliography{bibli}

\end{document}